\documentclass{article}
\usepackage[a4paper, margin=1in]{geometry}
\usepackage{graphicx} 
\usepackage{amsmath,amsfonts,amssymb,mathtools}
\usepackage{booktabs,makecell,multirow,array,colortbl,dcolumn,stfloats}
\usepackage[numbers]{natbib}
\usepackage{subcaption}
\usepackage{pgfplots}
\usepackage{tikz}
\usepackage{xcolor}
\usepackage{longtable}
\usetikzlibrary{arrows.meta, positioning}
\pgfplotsset{compat=1.18}
\usepackage{url} 
\usepackage[pdfusetitle,hidelinks]{hyperref} 
\usepackage[nameinlink]{cleveref} 
\usepackage{authblk}
\usepackage{orcidlink}

\Crefname{figure}{Fig.}{Figs.}
\Crefname{table}{Tab.}{Tabs.}
\Crefname{equation}{Eq.}{Eqs.}
\newtheorem{theorem}{Theorem}

\definecolor{cognac}{HTML}{D35400}
\definecolor{sage}{HTML}{45B39D}
\definecolor{dark_sage}{HTML}{117A65}
\definecolor{beige}{HTML}{FFA971}
\definecolor{golden_yellow}{HTML}{F3E4BC}

\title{\bf{Deep operator network models for predicting post-burn contraction}}
\author[1]{Selma Husanovic}
\author[2,3,4]{Ginger Egberts\orcidlink{0000-0003-3601-6496}}
\author[1]{Alexander Heinlein\orcidlink{0000-0003-1578-8104}\thanks{Corresponding author; E-mail: \texttt{a.heinlein@tudelft.nl}}}
\author[2]{Fred Vermolen\orcidlink{0000-0003-2212-1711}}

\affil[1]{\small Delft Institute of Applied Mathematics (DIAM), Faculty of Electrical Engineering, Mathematics \& Computer Science, Delft University of Technology, Delft, The Netherlands}
\affil[2]{\small Research Group Computational Mathematics (CMAT), Department of Mathematics and Statistics, University of Hasselt, Hasselt, Belgium}
\affil[3]{\small Department of Plastic, Reconstructive \& Hand Surgery, Amsterdam UMC, Amsterdam, The Netherlands}
\affil[4]{\small Dutch Burns Foundation, Beverwijk, The Netherlands}
\date{November 2024}

\begin{document}

\maketitle

\begin{abstract}
Burn injuries present a significant global health challenge. Among the most severe long-term consequences are \textit{contractures}, which can lead to functional impairments and disfigurement. Understanding and predicting the evolution of post-burn wounds is essential for developing effective treatment strategies. Traditional mathematical models, while accurate, are often computationally expensive and time-consuming, limiting their practical application. Recent advancements in machine learning, particularly in deep learning, offer promising alternatives for accelerating these predictions. This study explores the use of a \textit{deep operator network} (DeepONet), a type of neural operator, as a surrogate model for finite element simulations, aimed at predicting post-burn contraction across multiple wound shapes. A DeepONet was trained on three distinct initial wound shapes, with enhancement made to the architecture by incorporating initial wound shape information and applying sine augmentation to enforce boundary conditions. The performance of the trained DeepONet was evaluated on a test set including finite element simulations based on convex combinations of the three basic wound shapes. The model achieved an $R^2$ score of $0.99$, indicating strong predictive accuracy and generalization. Moreover, the model provided reliable predictions over an extended period of up to one year, with speedups of up to 128-fold on CPU and 235-fold on GPU, compared to the numerical model. These findings suggest that DeepONets can effectively serve as a surrogate for traditional finite element methods in simulating post-burn wound evolution, with potential applications in medical treatment planning. 
\end{abstract}

\medskip

\textbf{Keywords:} Operator learning, neural networks, surrogate model, wound modeling

\section{Introduction}
\noindent
Skin is the largest organ of the human body, and its weight is about 15 percent of the total body weight. Roughly speaking, skin consists of three layers: the epidermis (top layer), the dermis, and the subcutis. The subcutis is connected to underlying fat, organs, and muscles. If only the epidermis is damaged, then full healing will occur by the proliferation and migration of keratinocytes, which populate the epidermis. If deeper layers, such as the dermis or even the subcutis, are damaged, then the repair mechanism entails the regeneration of the dermis through collagen secretion by fibroblasts and the regeneration of the small blood vessel (capillary) network by revascularization. Serious burn injuries are characterized by damage on multiple skin (dermal) layers and over large areas of the human body. The healing of these large burn injuries often goes along with {\em dermal contraction} and/or {\em hypertrophy}. These mechanisms are mainly caused by the (strain-induced) differentiation of fibroblasts to myofibroblasts. While fibroblasts produce collagen type 1, the natural collagen in undamaged (embryonic) skin, myofibroblasts mainly secrete collagen type 3 at a large turnover per unit of time. Furthermore, myofibroblasts exert pulling forces, eventually making the wound area contract. At the same time, the epidermis is closed by the regeneration and migration of keratinocytes. Clinicians refer to the injury as a {\em wound} during epidermal healing, after which the injury is referred to as a {\em scar}. Although epidermal closure and dermal contraction influence each other, their exact interplay is still unclear, and therefore, we will not consider this interaction in the current study.

Dermal contraction may reduce the dermal pliability to such an extent that joints lose their mobility, and patients can no longer exercise or perform daily routines. In these cases, the occurred contraction is referred to as a {\em contracture}. Another side effect that patients with serious burn injuries endure is hypertrophy of scars. This hypertrophy is characterized by a changed dermal topography, which generally changes the appearance of the skin. 

While serious burn injuries often resulted in patient death in the past, nowadays, more and more patients survive severe burns. This is why current health care has shifted its focus from survival to the increase of quality of life of the survived patients. Therefore, clinical treatments are directed towards changing cellular behavior such that contraction and deposition of type 3 collagen are mitigated as much as possible. This mitigation is often achieved through methods such as dressings, clothing, skin grafting (skin transplantation by depositing epidermal parts from other body parts on the injured part), and splinting.

Many scientific studies, both clinical and laboratory, have been carried out and are ongoing to improve health care for burn injuries. In order to improve health care, one aims to steer the behavior of involved cells, and therefore, it is crucially important to understand the chain of biophysical mechanisms responsible for the evolution of post-burned skin. This knowledge needs to be validated using laboratory experiments and clinical observations, which represent data in numbers and patterns, and hence, the developed theory needs to be quantified. This quantification proceeds by using mathematical models based on either stochastic (random) processes and partial differential equations (continuum models) or by combinations of the two. To this extent, the current paper involves mathematical models in terms of partial differential equations (PDEs). These equations contain input parameters, such as elasticity and cell division rates, that are unknown on forehand due to patient-to-patient variations. For this reason, one simulation run cannot be deemed reliable in general, and therefore, one merely aims at predicting probabilities of clinical successes depending on the therapy. This can only be done by sampling the input parameters from prior probability distributions so that posterior probability distributions for the output parameters, such as the extent of wound contraction or strain energy, can be estimated. 

One of the bottlenecks in making current simulation tools applicable to clinical practice is the relatively long simulation time required to run one sample. To address this, neural network-based frameworks are being developed, where a regression-like model links input parameters to the output parameters of interest. Ideally, the resulting model can generate simulation runs at very large speedups compared to classical numerical simulations. 
The first neural networks to reproduce post-burn contraction finite element simulations in one and two dimensions have been introduced in~\cite{Egberts2022} and~\cite{Egberts2023}. These studies indicate high fidelity and high-speed reproduction of the finite element simulations. However, a limitation of these works is that all simulations were performed on a fixed spatial domain, meaning the wound's size and shape remain constant. This is not very representative of real-world scenarios. 

In this study, we propose using \textit{neural operators} to generate a surrogate model that reproduces the results of finite-element simulations, predicting post-burn wound evolution over time. Neural operators are a class of deep learning architectures with the primary application of learning surrogate maps for the solution operators of PDEs. Specifically, we will focus on \textit{deep operator networks (DeepONets)} for learning operators accurately and efficiently from a relatively small dataset~\cite{lu_deeponet_2021}. In particular, we aim to learn the solution operator of the two-dimensional \textit{morphoelastic model}, which describes post-burn skin evolution, thus accurately reproducing the finite element predictions. A key aspect of this research is incorporating multiple initial wound shapes. We train the neural operator network on a set of basic wound shapes, and evaluate its fidelity using convex combinations of the basic wound shapes as input. In this way, our research aims to extend the works \cite{Egberts2022} and \cite{Egberts2023}.

The paper is organized as follows. \Cref{sec:methodology} describes the methodology, focusing on numerical simulations and the machine learning models used. \Cref{sec:datasets_training} introduces the datasets and training setup for our neural network. The results are presented in \Cref{sec:results}. Finally, \Cref{sec:conclusion} provides the conclusions and discussion.     

\section{Methodology} \label{sec:methodology}
\noindent
This sections details the methodology, describing both the numerical simulations and the machine learning models employed. 

\subsection{Numerical simulations} 
\label{subsec:num_sim} \noindent
We first present the mathematical model for post-burn contraction, starting with the general system of equations and then providing a
more detailed description of the relevant biological and mechanical components. Furthermore, we give the boundary conditions, initial conditions, and numerical solver, and introduce an important measure of contraction. 

\paragraph{Mathematical morphoelastic model}
The general morphoelastic model for post-burn contraction was developed by Koppenol \cite{koppenol2017mathematical}, who used the theory of morphoelasticity developed by Hall \cite{hall2008modelling} to incorporate the formation of long-term deformations (contraction) into the dermal layer of the skin.
\par The model considers four biological constituents and three mechanical components as the primary variables. The biological constituents are the fibroblasts $(N)$, the myofibroblasts $(M)$, a generic signaling molecule $(c)$, and collagen $(\rho)$. The mechanical components are the dermal layer displacement $(\mathbf{u})$, the dermal layer displacement velocity  $(\mathbf{v})$, and the effective strain $(\pmb{\varepsilon})$. The following system of partial differential equations is used as a basis for the model: 
\begin{align}
    \frac{\mathrm{D} z_i}{\mathrm{D} t} + z_i (\nabla \cdot \mathbf{v}) &= - \nabla \cdot \mathbf{J}_i + R_i \label{eq:conservation_density}, \\
    \rho_t \left(\frac{\mathrm{D} \mathbf{v}}{\mathrm{D} t} + \mathbf{v}(\nabla \cdot \mathbf{v})\right) &= \nabla \cdot \pmb{\sigma} + \mathbf{f} \label{eq:conservation_momentum}, \\
    \frac{\mathrm{D} \pmb{\varepsilon}}{\mathrm{D} t} + \pmb{\varepsilon} \, \textrm{skw}(\nabla \mathbf{v}) - \textrm{skw}(\mathbf{\nabla \mathbf{v}}) \, \pmb{\varepsilon} \nonumber \\  + (\mathrm{tr}\left(\pmb{\varepsilon})-1\right)\mathrm{sym}(\mathbf{\nabla \mathbf{v}}) &= - \mathbf{G}. \label{eq:morphoelasticity}    
\end{align}
The operator $\frac{\mathrm{D}}{\mathrm{D} t}$ stands for the \textit{material derivative}:
\begin{equation*}
    \frac{\mathrm{D}}{\mathrm{D} t} \equiv \frac{\partial}{\partial t} + \mathbf{v}\cdot \nabla. 
\end{equation*} 
\Cref{eq:conservation_density} is the conservation equation for the cell density/concentration for each of the four biological constituents. Here, $z_i$ represents the concentration, $\mathbf{J}_i$ is the flux per unit area, and $R_i$ is a reaction term representing the kinetics of constituent $i$, for $i \in \{N, M, c, \rho\}$. 
\par \Cref{eq:conservation_momentum} is the conservation equation for linear momentum, where $\rho_t$ represents the total mass density of the dermal tissue, $\pmb{\sigma}$ is the stress tensor, and $\mathbf{f}$ is the total body force working on the dermal layer. Furthermore, $\mathbf{v} = \frac{\mathrm{D}\mathbf{u}}{\mathrm{D}t}$. We note that \Cref{eq:conservation_momentum} actually gives rise to multiple equations, one for each component of the velocity vector $\mathbf{v}$. For example, in 2D, this results in two equations.  
\par Finally, \Cref{eq:morphoelasticity} is the evolution equation
that describes how the infinitesimal effective strain ($\pmb{\varepsilon}$) changes
over time. This equation captures the morphoelasticity of the dermal layer, taking into account permanent deformation (in this case, contraction) and residual stresses. It was formulated by Hall \cite{hall2008modelling} and is based on his extensive theory on the zero stress state and morphoelasticity. The second-order tensor $\mathbf{G}$ is a growth tensor that describes the rate of
active change of the effective strain.

\paragraph{The fibroblasts}
To simplify notation, we will from now on replace $z_i$ by $i$. Hence, $z_N$ becomes $N$, the cell density of the fibroblasts in the dermis. The appropriate flux-term $\mathbf{J}_N$ incorporates both the random movement of fibroblasts through the dermal layer and the directed movement of fibroblasts up the gradient of signaling molecule $c$, if present. The former is modelled by a cell density-dependent Fickian diffusion, and the latter process is modelled using a simple model for chemotaxis \cite{hillen2009user}. Taken together, this gives
\begin{equation*}
    \mathbf{J}_N = -D_F F \nabla N + \chi_F N \nabla c, 
\end{equation*}
where $F = N + M$, $D_F$ is the (myo)fibroblast diffusion parameter, and $\chi_F$ is the chemotactic parameter. 
\par \Cref{eq:conservation_density} also contains a reaction term $R_N$ describing the kinetics of the fibroblasts. Three processes are considered: proliferation, differentiation into myofibroblasts, and apoptosis. The first is modelled using an adjusted logistic growth model. The presence of a signaling molecule $c$ is assumed to enhance both proliferation and cell differentiation:
\begin{equation*} 
    R_N = r_F\left(1 + \frac{r^{\mathrm{max}}_Fc}{a^I_c + c}\right)(1-\kappa_F F)N^{1+q} - k_F c N - \delta_N N. 
\end{equation*}
Here, the parameter $r_F$ is the cell division rate, $r^{\mathrm{max}}_F$ is the maximum factor with which the cell division rate can be enhanced due to the presence of the signaling molecule, and $a^I_c$ is the concentration of the signaling molecule that causes the half-maximum enhancement of the cell division rate. Furthermore, $\kappa_F F$ represents the reduction in the cell division rate due to crowding, $q$ is a fixed constant, $k_F$ is the signaling molecule-dependent cell differentiation rate of fibroblasts into myofibroblasts, and $\delta_N$ is the apoptosis rate of fibroblasts. 

\paragraph{The myofibroblasts}
For the myofibroblasts, the flux term in \Cref{eq:conservation_density} is very similar to the one for the fibroblasts:
\begin{equation*} 
    \mathbf{J}_M = -D_F F \nabla M + \chi_F M \nabla c. 
\end{equation*}
\par The reaction term describing the kinetics of myofibroblasts is also very similar. Almost the same adjusted logistic growth model is used, the only difference being the assumption that myofibroblasts solely divide when the generic signaling molecule is present:
\begin{equation*}
    R_M = r_F\left(\frac{(1+r^{\mathrm{max}}_F)c}{a^I_c + c}\right)(1-\kappa_F F)M^{1+q} - k_F c M - \delta_M M,     
\end{equation*}
where $\delta_M$ is the apoptosis rate of myofibroblasts. 

\paragraph{The signaling molecules}
We assume that the signaling molecules diffuse through the dermis according to linear Fickian diffusion:
\begin{equation*} 
    \mathbf{J}_c = -D_c \nabla c,
\end{equation*}
where $D_c$ is the diffusion coefficient of the generic signaling molecule. 
\par Furthermore, we assume that both fibroblasts and myofibroblasts release and consume the signaling molecules. Additionally, signaling molecules are removed from the dermis through proteolytic breakdown (breakdown of proteins into smaller components). Hence, we obtain the reaction term
\begin{equation*} 
    R_c = \frac{k_c(N+\eta^I M)c}{a^{II}_c + c} - \delta_c \, g(N,M,c,\rho)c,
\end{equation*}
where $k_c$ is the maximum net secretion rate of the signaling molecule, $\eta^I$ is the ratio of myofibroblasts to fibroblasts in the maximum net secretion rate of the signaling molecules and the collagen molecules, $a^{II}_c$ is the concentration of the signaling molecule that causes the half-maximum net secretion rate of the signaling molecule, and $\delta_c$ is the proteolytic breakdown rate of the signaling molecule. 
\par The function $g(N,M,c,\rho)$ represents the equilibrium concentration of a generic metalloproteinase (MMP). This enzyme is assumed to remove the signaling molecules through a proteolytic breakdown. In this study, we take the following relationship: 
\begin{equation*} 
    g(N,M,c,\rho) = \frac{(N+\eta^{II}M)\rho}{1+a_c^{III}c}.
\end{equation*}
The parameter $\eta^{II}$ is the ratio of myofibroblasts to fibroblasts in the secretion rate of the MMPs, and the $1/(1+a_c^{III}c)$ term represents the inhibition of the secretion of the MMPs due to the presence of the signaling molecule.   

\paragraph{The collagen molecules}
For collagen, we assume no active transport in the dermis. This means that the flux-term in \Cref{eq:conservation_density} is zero:
\begin{equation*} 
    \mathbf{J}_\rho = \mathbf{0}. 
\end{equation*}
The reaction term includes three processes: collagen molecules are produced by both fibroblasts and myofibroblasts, the secretion rate is enhanced in the presence of the signaling molecule, and there is an MMP-induced proteolytic collagen breakdown analogous to the removal of the signaling molecule:
\begin{equation*} 
    R_\rho = k_\rho \left(1 + \frac{k^{\mathrm{max}}_\rho c}{a^{IV}_c + c}\right)(N+ \eta^I
    M) - \delta_\rho \, g(N, M, c, \rho)\rho. 
\end{equation*}
Here, $k_\rho$ is the collagen molecule secretion rate, $k_\rho^{\mathrm{max}}$ is the maximum factor with which the secretion rate can be enhanced due to the presence of the signaling molecule, $a^{IV}_c$ is the concentration of the signaling molecule that causes the half-maximum enhancement of the secretion rate, and $\delta_\rho$ is the degradation rate of the collagen molecules.

\paragraph{The mechanical components}
In \Cref{eq:conservation_momentum}, a visco-elastic constitutive relation is used for the stress-strain relation in the dermal layer. The visco-elastic
relation for the dermal stress is: 
\begin{equation*}
    \pmb{\sigma} = \mu_1 \mathrm{sym}(\nabla \mathbf{v}) + \mu_2[\mathrm{tr}(\mathrm{sym}(\nabla \mathbf{v}))\mathbf{I}]\, + \frac{E \sqrt{\rho}}{1+\nu} \left(\pmb{\varepsilon} + \mathrm{tr}(\pmb{\varepsilon})\frac{\nu}{1-2\nu}\mathbf{I}\right),
\end{equation*} 
where $\mu_1$ and $\mu_2$ are the shear and bulk viscosity, respectively, and $\nu$ is the Poisson's ratio. Additionally, $E\sqrt{\rho}$ represents Young’s modulus (stiffness), which we assume to be dependent on the concentration of the collagen molecules.  
\par Additionally, the total body force $\mathbf{f}$ in \Cref{eq:conservation_momentum} needs a more precise description. We assume that the myofibroblasts generate an isotropic stress due to their pulling on the extracellular matrix, which is proportional to the product of the cell density of the myofibroblasts and a simple function of the concentration of the collagen molecules:
\begin{align*} 
    \mathbf{f} &= \nabla \cdot \pmb{\psi}, \\
    \pmb{\psi} &= \xi M \left(\frac{\rho}{R^2 + \rho^2}\right) \mathbf{I}.
\end{align*}
Here, $\pmb{\psi}$ is a second-order tensor representing the total generated stress by the myofibroblast population, the parameter $\xi$ is the generated stress per unit cell density and the inverse of the unit collagen concentration, and $R$ is a fixed constant.  
\par Finally, we consider the growth contribution tensor $\mathbf{G}$ in \Cref{eq:morphoelasticity}. We assume that the rate of active change of the effective strain is proportional to the product of the amount of effective strain, the local concentration of the MMPs, the local concentration of the signaling molecule, and the inverse of the local concentration of the collagen molecules. Taken collectively, this results in the following symmetric tensor:
\begin{equation*} 
    \mathbf{G} = \zeta \left(\frac{g(N,M,c,\rho)c}{\rho}\right)\pmb{\varepsilon} = \zeta \left( \frac{(N+\eta^{II}M)c}{1+a_c^{III}c} \right) \pmb{\varepsilon},
\end{equation*}
where the parameter $\zeta$ is the rate of morphoelastic change.  

\paragraph{Boundary conditions} 
We locate the $xy$-plane parallel to the surface of the skin and write
\begin{equation*}
\boldsymbol{v}=\begin{bmatrix}v_1\\v_2\end{bmatrix},\qquad\text{and}\qquad \boldsymbol{\varepsilon}=\begin{bmatrix}\varepsilon_{11}&\varepsilon_{12}\\\varepsilon_{21}&\varepsilon_{22}\end{bmatrix}.
\end{equation*}
We perform the computations on an arbitrary skin depth. Such a configuration can approximate the kinetics of a wound on a non-curved body part. Hence, we disregard all dependencies on the depth of the burn into the skin.

We perform computations on a reduced, symmetrical domain to reduce the computational workload. The solution inherits this symmetrical property. We define the \emph{computational domain} by $\Omega_{\mathbf{x}}$ and the \emph{boundary of the computational domain} by $\partial\Omega_{\mathbf{x}}$. The distance between $\partial\Omega_{\mathbf{x}}$ and the modelled burn is sufficiently large to prevent too much variable diffusion near the boundary. Furthermore, the computational domain implicitly depends on the time $t$ given that $\mathbf{x}=\mathbf{x}(t)$. We do not specify any boundary conditions for $\rho$ and $\boldsymbol{\varepsilon}$ because of over-determination, since the equations for $\rho$ and $\boldsymbol{\varepsilon}$ are ordinary differential equations for time $t$.

Let $\partial\Omega_{\mathbf{x}} = \{\partial\Omega_{\mathbf{x}}^o ,\partial\Omega_{\mathbf{x}}^h,\partial\Omega_{\mathbf{x}}^v$\}, where $\partial\Omega_{\mathbf{x}}^o$ represents the outer non-symmetrical boundaries, $\partial\Omega_{\mathbf{x}}^h$ represents the horizontal symmetrical boundary where $y=0$, and $\partial\Omega_{\mathbf{x}}^v$ represents the vertical symmetrical boundary where $x=0$. For the chemicals, the following boundary conditions hold for all time $t$ and all
\begin{align*}
\mathbf{x}\in\partial\Omega_{\mathbf{x}}^o&: \quad N(\mathbf{x};t)=\overline{N},\quad M(\mathbf{x};t)=\overline{M},\\
&\quad\quad\quad\text{and}\quad c(\mathbf{x};t)=\overline{c},\\
\mathbf{x}\in\partial\Omega_{\mathbf{x}}^p&: \quad \mathbf{J}_{N/M/c} \cdot \mathbf{n} = 0,
\end{align*}
where $p\in\{h,v\}$ and $\mathbf{n}$ is the outward pointing normal vector. We use similar conditions for the mechanics, for all time $t$ and all
\begin{align*}
\mathbf{x}\in\partial\Omega_{\mathbf{x}}^o&: \quad \boldsymbol{v}(\mathbf{x};t)=0,\\
\mathbf{x}\in\partial\Omega_{\mathbf{x}}^p&: \quad \boldsymbol{v}\cdot \mathbf{n} = 0 \quad \text{and} \quad (\sigma \cdot \mathbf{n})\cdot\tau = 0,
\end{align*}
where $\tau$ is the tangential vector.

\paragraph{Initial conditions}
The initial conditions describe the cell densities and the concentrations at the onset of the proliferative phase of wound healing. Signaling molecules are present in the wound because of their secretion in the  inflammatory phase of wound healing. Furthermore, fibroblasts and collagen are initially assumed to be present in the wound, whereas myofibroblasts are assumed to be absent. 
\par The initial wounded area is denoted by $\Omega^w(0) \subset \Omega_{\mathbf{x},0}$. The unwounded area is then $\Omega_{\mathbf{x},0}\setminus\Omega^w(0)$. Let $d(\mathbf{x})$ be the shortest distance from a point $\mathbf{x}\in\Omega^w$ to the wound boundary. Let $\Omega^w_s = \{{\bf x} \in \Omega^w(0)~:~d({\bf x}) \ge s\}$. Then, for $z\in\{N,c,\rho\}$ we have the following initial densities or concentrations: 
\begin{align*}
    z(\mathbf{x},0) = 
    \begin{cases}
        \tilde{z}, \quad &\mathbf{x}\in \Omega_s^w,\\
        \overline{z}, \quad &\mathbf{x}\in \Omega_{\mathbf{x},0}\setminus\Omega^w(0),
    \end{cases}
\end{align*}
where $\tilde{z}, \overline{z} \in \mathbb{R}^+$, the latter indicating the equilibrium value. Furthermore, 
\begin{equation*}
    M(\mathbf{x},t) = \overline{M} = 0, \quad \mathbf{x}\in\Omega_{\mathbf{x},0}.
\end{equation*}
For the wound boundary steepness, we use half a period of sine-functions for $N$, $c$, and $\rho$ to smoothly transition from the wound to the unwounded area. 
 
Regarding the initial conditions for the mechanical part of the model, all quantities are initialized with zero.

\paragraph{Numerical discretizations and solvers}
We solve the model equations by applying the finite element method. We note that we add the term $\varepsilon_{i,j}[\nabla\cdot\boldsymbol{v}]$ for $i,j\in\{1,2\}$ to the left-hand side and the right-hand side of the effective Eulerian strain equations before deriving the weak formulation. Furthermore, we make use of $\varepsilon_{21}=\varepsilon_{12}$ (symmetric effective strain tensor), which we have proved for all time $t$ \cite{Egberts2020a}. For a complete overview of the derivation of the finite element formulations, we refer to the appendix in Koppenol’s PhD thesis \cite{KoppenolThesis}. 

We subdivide the computational domain into a finite number of non-overlapping triangles $\Delta_p$ (i.e., the elements) that are as equilateral as possible. Let $X_h(t) \subset H^1(\Omega_{\mathbf{x}})$ the finite element subspace and $\mathbf{a}_j,j\in\{1,\dots,n\}$, $n\in\mathbb{N}$ the coordinates of the vertices of the elements. We choose piecewise linear Lagrangian basis functions $\varphi_i\in X_h(t)$ with $\varphi_i(\mathbf{a}_j;t)=\delta_{ij}$ for $i,j\in\{1,\dots,n\}$ as basis functions for the finite-dimensional subspace $X_h(t)$, where $\delta_{ij}$ denotes the Kronecker delta function.

We approximate the local displacements of the dermal layer $(\boldsymbol{u})$ with
\begin{equation*}
{\boldsymbol{u}}({\bf x};t+\Delta t) \approx {\boldsymbol{u}}({\bf x};t) + {\boldsymbol{v}}({\bf x};t)\Delta t.
\end{equation*}
We use the initial condition $\boldsymbol{u}(\mathbf{x};0) = 0,\quad\forall\mathbf{x}\in\Omega_{\mathbf{x},0}$. 
Furthermore, we update the mesh (triangulation) in every time integration step, and we determine the quality of this updated mesh by computing
$$\min_{e_k}\left|\mathbf{J}_{e_k}\right|/\max_{e_k}\left|\mathbf{J}_{e_k}\right|,\quad e_k\in\Omega,$$ with $\mathbf{J}$ the Jacobian. In case $\min_{e_k}\left|\mathbf{J}_{e_k}\right|/\max_{e_k}\left|\mathbf{J}_{e_k}\right|<0.5$, we perform (global) remeshing. Note that the Jacobian is twice the area of a triangular element. 

We use mass lumping and a semi-implicit flux corrected transport limiter \cite{Moeller2008} that enforces the positiveness of solutions, so that loss of monotonicity (i.e., spurious oscillations) is suppressed.

\paragraph{Relative surface area} During wound healing, due to myofibroblasts pulling on the surrounding collagen fibers, the wound contracts towards its center and retracts after these cells disappear. \textit{The relative surface area} (of the) \textit{wound} (RSAW) is an important measure, as it gives valuable information about contraction. It is defined as follows:
\begin{equation*}
    RSAW(t) = \frac{area(\Omega^w(t))}{area(\Omega^w(0))}. 
\end{equation*}
The minimum RSAW value corresponds to maximum contraction, and the asymptotic value indicates the long-term (permanent) contraction.  

\subsection{Machine learning models}
\noindent
In this section, we introduce the methodology of our machine learning models for predicting the evolution of the wound shape over time. Therefore, we first provide a brief introduction to neural networks, then introduce the concept of operator learning, and finally, introduce the network architecture of our machine learning model, which is based on DeepONets~\cite{lu_deeponet_2021}. 

\paragraph{Neural networks}

Let us consider a generic supervised machine learning task. Therefore, let input data $X = (x_1,\ldots,x_N)$ and output data $Y=(y_1,\ldots,y_N)$ be given, with $x_i \in \mathbb{R}^n$ and $y_i \in  \mathbb{R}^m$, for $1 \leq i \leq N$. Then, we aim at finding a parametrized function 
$$
	f_\theta : \mathbb{R}^n \to \mathbb{R}^m,
$$
such that 
\begin{equation} \label{eq:fit}
	f_\theta (x_i) \approx y_i;
\end{equation}
where $\theta \in \mathbb{R}^K$ is a vector consisting of all parameters of the model. 

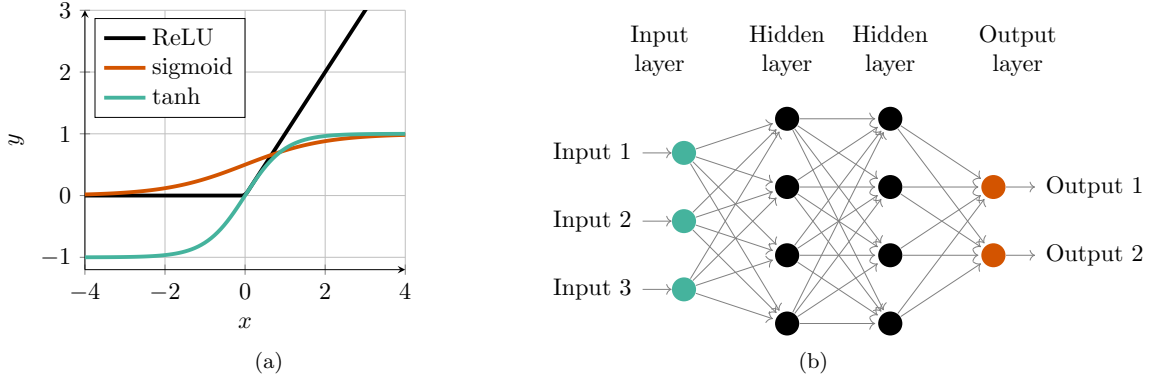
\begin{figure}[t]
    \centering
    \resizebox{0.9\textwidth}{!}{
    \begin{subfigure}[t]{0.49\textwidth}
        \begin{tikzpicture}
        \begin{axis}[
          width=0.8\textwidth,
          xlabel={$x$},
          ylabel={$y$},
          domain=-4:4,
          samples=100,
          axis lines=left,
          grid=both,
          xmin=-4, xmax=4,
          ymin=-1.2, ymax=3,
          legend pos=north west,
          legend cell align=left,
        ]
        
        \addplot[black,ultra thick] {max(0,x)};
        \addlegendentry{ReLU}
        
        \addplot[cognac,ultra thick] {1/(1 + exp(-x))};
        \addlegendentry{sigmoid}
        
        \addplot[sage,ultra thick] {tanh(x)};
        \addlegendentry{tanh}
        
        \end{axis}
        \end{tikzpicture}
        \caption*{(a)}
    \end{subfigure}
    \begin{subfigure}[t]{0.49\textwidth}
        \centering
        \begin{tikzpicture}[shorten >=1pt,->,draw=black!50, node distance=1.5cm]
        \tikzstyle{neuron}=[circle,fill=black!25,minimum size=10pt,inner sep=0pt]
        \tikzstyle{input neuron}=[neuron, fill=sage];
        \tikzstyle{output neuron}=[neuron, fill=cognac];
        \tikzstyle{hidden neuron}=[neuron, fill=black];
        \tikzstyle{annot} = [text width=4em, text centered]
    
        \foreach \name / \y in {1,...,3}
          \node[input neuron, pin={[pin edge={<-}]left:Input \y}] (I-\name) at (0,-\y) {};
    
        \foreach \name / \y in {1,...,4}
          \node[hidden neuron] (H1-\name) at (1.5,-\y+0.5) {};
    
        \foreach \name / \y in {1,...,4}
          \node[hidden neuron] (H2-\name) at (3,-\y+0.5) {};
    
        \foreach \name / \y in {1,...,2}
          \node[output neuron, pin={[pin edge={->}]right:Output \y}] (O-\name) at (4.5,-\y-0.5) {};
    
        \foreach \source in {1,...,3}
          \foreach \dest in {1,...,4}
            \path (I-\source) edge (H1-\dest);
    
        \foreach \source in {1,...,4}
          \foreach \dest in {1,...,4}
            \path (H1-\source) edge (H2-\dest);
    
        \foreach \source in {1,...,4}
          \foreach \dest in {1,...,2}
            \path (H2-\source) edge (O-\dest);
    
        \node[annot, above of=H1-1, node distance=1cm] (hl1) {Hidden layer};
        \node[annot, above of=H2-1, node distance=1cm] (hl2) {Hidden layer};
        \node[annot, left=0.2cm of hl1] {Input layer};
        \node[annot, right=0.2cm of hl2] {Output layer};
        \end{tikzpicture}
        \caption*{(b)}
    \end{subfigure}}
    
    \caption{Nonlinear activation functions (a) and a feedforward neural network with two hidden layers (b)}
    \label{fig:sigma_nn}
\end{figure}

Neural networks are a specific class of parametrized functions. In particular, neural networks are composed of functions of the form
\begin{equation} \label{eq:fj}
	f_{W_j,b_j}^{(j)} (x) \coloneq \sigma ( W_j \cdot x + b_j )
\end{equation}
with the index $j$, $W_j \in \mathbb{R}^{n_j \times n_{j+1}}$ and $b_j \in \mathbb{R}^{n_{j+1}}$, and a nonlinear \emph{activation function} $\sigma$ which is applied component-wise to a vector; typical examples for the activation function $\sigma$ are shown in~\Cref{fig:sigma_nn} (a). Here, we consider only $\sigma(x) = \max(0,x)$, (ReLU). The function $f^{(j)}_{W_j,b_j}$ is a composition of a multivariate affine function and a nonpolynomial activation function $\sigma$. A simple feedforward neural network, or \emph{multi-layer perceptron (MLP)}, is then the composition of multiple functions of the form given in \Cref{eq:fj} and one final linear map, that is,
$$
	f_\theta \coloneq A \circ f_{W_l,b_l}^{(l)} \circ f_{W_{l-1},b_{l-1}}^{(l-1)} \circ \cdots \circ f_{W_{1},b_{1}}^{(1)}.
$$ 
Here, $A \in \mathbb{R}^{m \times n_{l+1}}$ corresponds to the linear \emph{output layer}, and each function $f_j$, for $1 \leq j \leq l$, corresponds to one so-called \emph{hidden layer} of the model; see also \Cref{fig:sigma_nn} (b). The matrices and vectors $W_j$ and $b_j$, for $1 \leq j \leq l$, as well as the matrix $A$ constitute the parameter vector $\theta$. We call a neural network \textit{deep}, if the number of hidden layers is sufficiently large, typically for at least $3$ hidden layers.

In order for the model $f_\theta$ to fit the data as in \Cref{eq:fit}, a minimization problem is formulated. In the simplest case of regression, fitting the model can be equivalently written as finding
$$
	\theta^\star = \arg\min_\theta \mathcal{L}(\theta)
	\quad \text{with} \quad
	\mathcal{L}(\theta) = \frac{1}{N} \sum_{i = 1}^N \| f_\theta (x_i) - y_i \|_2^2.
$$
We also call $\mathcal{L}(\theta)$ the \emph{loss function}, which here is based on the mean squared error (MSE). Due to the nonlinearity of $f_\theta$, the loss function is complicated and, in general, non-convex. The minimization is usually performed using a gradient-based optimizer or using a quasi-Newton method. Typical examples are stochastic gradient descent (SGD), and its variant Adam~\cite{kingma_adam_2017}, or the limited memory Broyden--Fletcher--Goldfarb--Shanno (L-BFGS) method~\cite{liu_limited_1989}. These methods require gradient information, which can be computed efficiently using automatic differentiation using the backpropagation algorithm~\cite{kelley_gradient_1960}. 

For more details, see for instance ~\cite{goodfellow_deep_2016}. Note that neural networks can also be used for semi-supervised and unsupervised learning. However, we will only consider the case of supervised learning in this paper.

\paragraph{Operator learning}
In~\cite{kovachki_neural_2022}, neural operators have been introduced  as neural networks for approximating operators, that is, maps between infinite-dimensional function spaces. Suitable neural operator architectures have universal approximation properties, meaning that, with sufficient model capacity (i.e., a sufficiently large number of model parameters), they can approximate any given nonlinear continuous operator up to arbitrary precision. Furthermore, in~\cite{kovachki_neural_2022}, neural operators are discretization-invariant: 1) the model can act on any discretization of the input function; 2) the model can be evaluated at any point in the output domain; 3) when refining the discretization, the model converges to a continuum operator.

Neural operators are particularly well-suited for creating surrogate models for numerical solvers for PDEs, since the input of the model can be a function, such as
\begin{itemize}
	\item a right-hand side function of a PDE, 
	\item a spatial distribution of a material parameter, 
	\item a function describing a boundary condition, or 
	\item a level-set function describing the geometry,
\end{itemize}
and the output of the neural operator can be the solution of the initial boundary value problem (IBVP). This is different from, for instance, classical physics-informed neural networks (PINNs)~\cite{raissi_physics-informed_2019}, which generally only approximate the solution of a single IBVP. 

Recently, many different neural operator architectures have been introduced, with Fourier neural operators (FNOs)~\cite{li_fourier_2021} and DeepONets~\cite{lu_deeponet_2021} being the most popular choices. We also highlight convolutional neural operators (CNOs)~\cite{raonic_convolutional_2023}, which satisfy structure-preserving continuous-discrete equivalence enabling learning operators without discretization-dependent aliasing errors. The proposed CNO architecture is an extension of the popular U-Net architecture~\cite{10.1007/978-3-319-24574-4_28}. 

\begin{figure}[t]
    \centering
    \resizebox{0.5\textwidth}{!}{
    \begin{tikzpicture}[scale=0.5, >=Stealth]

    \node[draw, rectangle] (branch_inputs) {
        \begin{tabular}{c}
            $v(x_1)$ \\ 
            $v(x_2)$ \\ 
            \vdots \\
            $v(x_m)$
        \end{tabular}
    };

    \node[draw, rectangle, right=0.5cm of branch_inputs] (branch_net) {Branch Net};

    \node[draw, rectangle, right=0.5cm of branch_net] (branch_outputs) {
        \begin{tabular}{c}
            $b_1$ \\ 
            $b_2$ \\
            $\vdots$ \\ 
            $b_p$ \\ 
        \end{tabular}
    };

    \node[draw, rectangle, below=1cm of branch_inputs] (trunk_inputs) {
        \begin{tabular}{c}
            $y$
        \end{tabular}
    };

    \node[draw, rectangle, right=0.9cm of trunk_inputs] (trunk_net) {Trunk Net};  

    \node[draw, rectangle, right=0.65cm of trunk_net] (trunk_outputs) {
        \begin{tabular}{c}
            $t_1$ \\ 
            $t_2$ \\
            $\vdots$ \\ 
            $t_p$
        \end{tabular}
    };

    \node[circle, draw, right=0.6cm of branch_outputs, yshift=-1.2cm] (inner_product1) {$\bigotimes$}; 

    \node[draw, rectangle, right=0.5cm of inner_product1] (output1) {$G(v)(y)$};
    
    \draw[->] (branch_inputs) -- (branch_net);
    \draw[->] (branch_net) -- (branch_outputs);

    \draw[->] (trunk_inputs) -- (trunk_net);
    \draw[->] (trunk_net) -- (trunk_outputs);

    \draw[->] ([yshift=1cm]branch_outputs.east) -- (inner_product1.west);

    \draw[->] (trunk_outputs.east) -- (inner_product1.west);

    \draw[->] (inner_product1) -- (output1);

    \end{tikzpicture} }
	\caption{(Unstacked) DeepONet architecture; based on~\cite{lu_deeponet_2021}
	\label{fig:deeponet}
	}
\end{figure}
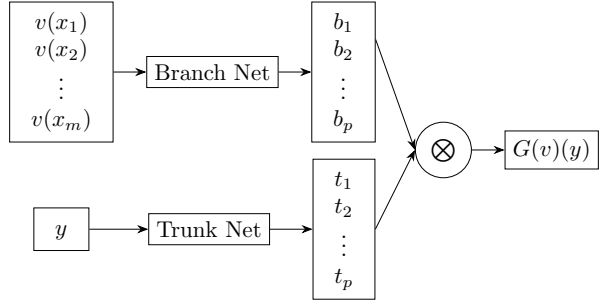
Here, we build a neural operator learning model based on the DeepONet architecture, and for the sake of brevity, we will focus only on a brief description of this approach. 
To this extent, consider a generic nonlinear continuous operator $G: V \to W$, where $V$ and $W$ are compact, nonempty subsets of Banach spaces, which amount to normed function spaces on a nonempty, compact subset of physical space-time. 
We are interested in evaluating $G(v)$ at discrete points in the nonempty compact portion of physical space-time. The user has to input two sets: (1) the points in a compact region in space-time where (s)he wants to evaluate $G(v)$ and (2) the input function at points in a compact region in space-time, that is $v(x_1), \dots, v(x_m)$. Note that the compact regions where $v$ is obtained and where $G(v)$ is determined do not have to be the same. The Universal Approximation Theorem for Operator Learning states that each nonlinear mapping between nonempty subsets of Banach spaces can be approximated arbitrarily well, provided that sufficient input information (such as the dimensionality and features of the input data) and complexity (the depth, width, and structure of the neural network) are incorporated:

\begin{theorem}[Universal approximation for operators~\cite{chen_universal_1995,lu_deeponet_2021}] \label{thm:approx}
	Let $\sigma$ be a continuous nonpolynomial function, $X$ be a Banach space, $K_1 \subset X$ and $K_2 \subset \mathbb{R}^d$ two compact subsets, respectively, $V$ be a compact subset of $\mathcal{C}(K_1)$, and $G$ be a nonlinear continuous operator $V \to \mathcal{C}(K_2)$. 
	Then, for any $\varepsilon > 0$, there exist positive integers $n,p,m$ and constants $c_i^k, w_{ij}^k, b_i^k, \hat{b}_k \in \mathbb{R}, \hat{w}_k \in \mathbb{R}^d, x_j \in K_1$, for $i=1,\ldots,n$, $k=1,\ldots,p$, and $j=1,\ldots,m$, such that
	\begin{equation} \label{eq:deeponet_ua}
		\big| 
		G(v)(y)
		-
		\sum_{k = 1}^{p} 
		\underbrace{
			\sum_{i = 1}^{n}
			c_i^k \sigma \big( \sum_{j = 1}^{m} w_{ij}^k u(x_j) + b_i^k \big) 
		}_{\text{branch}}
		\underbrace{
			\sigma (\hat{w}_k y + \hat{b}_k) 
		}_{\text{trunk}}
		\big|  
		< 
		\varepsilon
	\end{equation}
	for all $v \in V$ and $y \in K_2$.
\end{theorem}
The idea of DeepONets is then to replace what is denoted as \emph{branch} and \emph{trunk} in~\Cref{eq:deeponet_ua} by (deep) neural networks. Note that~\Cref{eq:deeponet_ua} relates to~\Cref{eq:fj} in the sense that $w_{ij}^k$ and $\hat{w}_k$ are the coefficients of weight matrices $W_j$ of the branch and trunk, and $b_i^k$ and $\hat{b}_k$ are the coefficients of bias vectors $b_j$. \Cref{thm:approx} warrants the existence of such a network; however, it says nothing about the rate of convergence. 

In its general form, the DeepONet architecture is depicted schematically in~\Cref{fig:deeponet}. Intuitively, the trunk network corresponds to spatial basis functions $t_1,\ldots,t_p$, whereas the branch network provides corresponding coefficients $b_1,\ldots, b_p$ in the basis representation 
$$
	G(v) (y) \approx \sum_{k=1}^{p} b_k (x_1,\ldots,x_m) t_k (y).
$$

Many variants of DeepONets have been introduced, where specific types of networks are chosen as the branch and trunk networks, for instance, U-Net~\cite{diab_u-deeponet_2023} and graph neural network (GNN)~\cite{sun_deepgraphonet_2022} architectures. In many related works, the basis functions are prescribed, instead of learning them via a neural network, for instance, using classical discretizations~\cite{wang_deep_2020,moller2021physics,meethal_finite_2023,griese_fem-based_2024}. Furthermore, the basis property has also been recently used to learn coarse basis functions for robust preconditioning~\cite{kopanicakova_deeponet_2024}. For improved precision, DeepONets can be stacked~\cite{howard_stacked_2023} and combined with domain decomposition techniques~\cite{heinlein_multifidelity_2024}.

\paragraph{Deep operator network architecture}
We propose a machine learning model, based on the DeepONet architecture, for predicting post-burn contraction over time. The exact architecture of our model is schematically depicted in \Cref{fig:architecture}. It consists of a branch and trunk network, both of which are simple MLPs with three hidden layers, containing 50 neurons per layer. 
The branch takes as input five parameters from the morphoelastic model describing wound contraction: the cell diffusion constant $D_F$, the chemotactic constant $\chi_F$, the signaling diffusion constant $D_c$, the cell differentiation rate $k_F$, and the half-maximum cell division enhancement rate $a_c^I$. These parameters describe the behavior of cells involved in the wound healing process and are characteristics that vary between patients. The branch network outputs a vector $(b_1, \dots, b_{2p}) \in \mathbb{R}^{2p}$.
\par
\begin{figure}[t]
    \centering
    \begin{subfigure}[t]{.32\textwidth}
        \centering
        \begin{tikzpicture}        
            \draw[->] (0,0) -- (4,0) node[anchor=north] {$x$};
            \draw[->] (0,0) -- (0,3) node[anchor=east] {$y$};

            \draw[thick] (3,0) -- (3,2) -- (0,2);
            \node[cognac, anchor=south west] at (3,2) {$(x_l, y_l)$};
            \fill[cognac] (3,2) circle (2pt);

            \draw[thick, sage] (2,0) node[below] {$x_{\text{cut}}$} -- (2,1) -- (0,1) node[left] {$y_{\text{cut}}$};
            \node[sage, anchor=south] at (2,1) {$(x_m, y_m)$};
            \fill[sage] (2,1) circle (2pt);
        \end{tikzpicture}
        \caption*{(a)}
    \end{subfigure}%
    \hfill
    \begin{subfigure}[t]{.32\textwidth}
        \centering
        \begin{tikzpicture}
            \draw[->] (0,0) -- (4,0) node[anchor=north] {$x$};
            \draw[->] (0,0) -- (0,3) node[anchor=east] {$y$};

            \draw[thick] (3,0) -- (3,2.25) -- (0,2.25);
            \node[cognac, anchor=south west] at (3,2.25) {$(x_l, y_l)$};
            \fill[cognac] (3,2.25) circle (2pt);

            \draw[thick, sage] (2,0) node[below] {$x_{\text{cut}}$} -- (0,1.5) node[left] {$y_{\text{cut}}$};
            \node[sage, anchor=south west] at (0.8,0.9) {$(x_m, y_m)$};
            \fill[sage] (0.8,0.9) circle (2pt);
        \end{tikzpicture}
        \caption*{(b)}
    \end{subfigure}%
    \hfill
    \begin{subfigure}[t]{.32\textwidth}
        \centering
        \begin{tikzpicture}
            \draw[->] (0,0) -- (4,0) node[anchor=north] {$x$};
            \draw[->] (0,0) -- (0,3) node[anchor=east] {$y$};

            \draw[thick] (3.75,0) -- (3.75,2.25) -- (0,2.25);
            \node[cognac, anchor=south west] at (3.75,2.25) {$(x_l, y_l)$};
            \fill[cognac] (3.75,2.25) circle (2pt);

            \draw[thick, sage] (2.5,0) node[below] {$x_{\text{cut}}$} arc[start angle=0, end angle=90, x radius=2.5, y radius=1.5] node[left] {$y_{\text{cut}}$};
            \node[sage, anchor=south west] at (1,1.375) {$(x_m, y_m)$};
            \fill[sage] (1,1.375) circle (2pt);        
        \end{tikzpicture}
        \caption*{(c)}
    \end{subfigure}
    \caption{The three initial wound shapes used for training: rectangle (a), rhombus (b), and ellipse (c). Due to symmetry, a quarter of the complete domain is considered. Additionally, a visualization of the quadruple $(y_{cut}, x_m, y_m, x_{cut})$, which serves as extra input to the trunk. The coordinate $(x_l, y_l)$ is required for the sine augmentation step}
    \label{fig:domain_and_wound}
\end{figure}
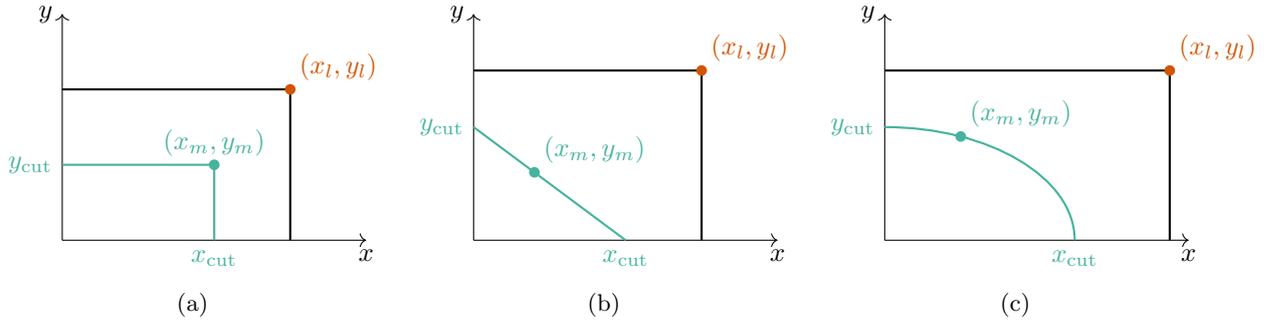

 The trunk network takes as input a coordinate $(t, x, y)$ in which we want to evaluate the displacement $\mathbf{u}(t,x,y) = (u_1, u_2)$. The trunk additionally takes as input information about the initial wound shape, given by the quadruple $(y_{cut},x_m,y_m,x_{cut})$. \Cref{fig:domain_and_wound} shows the three wound shapes used for training (see also \Cref{subsec:training}) and visualizes this extra input to the trunk. The latter uniquely determines the initial wound shape at $t=0$. The output of the trunk network is a vector $(c_1, \dots, c_{p}) \in \mathbb{R}^{p}$. In our simulations, we fix $p=50$. 
\par The output of the DeepONet is an inner product of the outputs of the trunk and branch networks. Specifically, let $T(\cdot)$ be the operation of the trunk and $B(\cdot)$ be the operation of the branch network, then we have
\begin{align*}
    c_1, \dots c_p &= T(t, x, y, y_{cut}, x_m, y_m, x_{cut}), \\
    b_1, \dots b_{2p} &= B(D_F,\chi_F, D_C,k_F,a_C^I), \\
    \hat{u}_1 &= \sum_{i=1}^p b_ic_i, \\
    \hat{u}_2 &= \sum_{i=1}^p b_{i+p}c_i. 
\end{align*}
\par
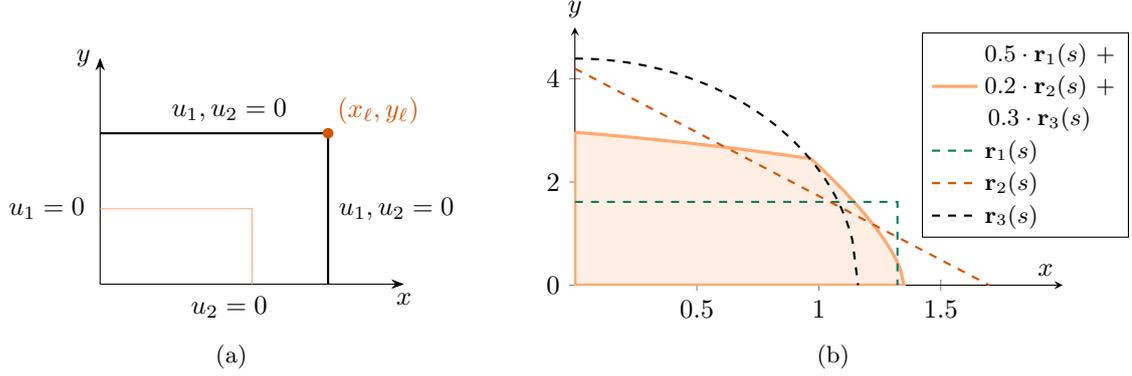
\begin{figure}[t]
\begin{subfigure}[t]{0.49\textwidth}
    \centering
        \begin{tikzpicture}[>=Stealth]       
            \draw[->] (0,0) -- (4,0) node[anchor=north] {$x$};
            \draw[->] (0,0) -- (0,3) node[anchor=east] {$y$};

            \draw[thick] (3,0) -- (3,2) -- (0,2);
            \node[cognac, anchor=south west] at (3,2) {$(x_{\ell}, y_{\ell})$};
            \fill[cognac] (3,2) circle (2pt);

            \draw[ultra thin, beige] (2,0) -- (2,1) -- (0,1); 

            \node at (1.7,2.3) {$u_1,u_2=0$};   
            \node at (3.9,1) {$u_1,u_2=0$};  
            \node at (-0.7,1) {$u_1=0$};  
            \node at (1.7,-0.3) {$u_2=0$};  
        \end{tikzpicture}
    \caption*{(a)}
\end{subfigure}
\begin{subfigure}[t]{0.49\textwidth}
    \centering
    \begin{tikzpicture}
      \begin{axis}[
          axis lines = middle,
          xmin = 0, xmax=2,
          ymin = 0, ymax = 5,
          xlabel = $x$,
          ylabel = $y$,
          ylabel style={anchor=south},
          xtick = {0, 0.5, 1, 1.5},
          extra y ticks={0},
          extra y tick labels={$0$},
          width=8cm, height=5cm,
          legend cell align=left,
          legend style={
            font = \small,
            xshift= 1cm, 
            },
        ]

        \addplot[beige, very thick, fill=beige!20] table {data/convex.txt};
        \addplot[dark_sage, dashed, thick] table {data/square.txt};
        \addplot[cognac, dashed, thick] table {data/rhombus.txt};
        \addplot[dashed, thick] table {data/ellipse.txt};

        \legend{
          {\shortstack{$0.5 \cdot {\bf r}_1(s)$ + \\ $0.2 \cdot {\bf r}_2(s)$ + \\$0.3 \cdot {\bf r}_3(s) \,\,\, \, \,$}},
          {${\bf r}_1(s)$},
          {${\bf r}_2(s)$},
          {${\bf r}_3(s)$}
        }
    
      \end{axis}
    \end{tikzpicture}
    \caption*{(b)}
\end{subfigure}
\caption{Enforcement of boundary conditions as a result of sine augmentation (a) and a convex combination of the three basic wound shapes (b)}
\label{fig:sine_and_bnd}
\end{figure}

 Finally, we apply \textit{sine augmentation} to the DeepONet outputs, ensuring zero displacements at the top and right boundaries of our domain. Moreover, it ensures zero $x$-displacement on the left boundary and zero $y$-displacement on the bottom boundary. \Cref{fig:sine_and_bnd} (a) summarises the effect on all four boundaries. This hard enforcement of boundary conditions is not new, examples can already be found in early works employing neural networks to solve PDEs; see, for instance,~\cite{lagaris_artificial_1998}. The final outputs of our model are given by
\begin{equation} \label{eq:sine}
	\begin{aligned}
		u_1 &= \hat{u}_1 * \sin\left(\frac{\pi}{x_l}x\right) \cos\left(\frac{\pi}{2y_l}y\right), \\
		u_2 &= \hat{u}_2 * \sin\left(\frac{\pi}{y_l}y\right) \cos\left(\frac{\pi}{2x_l}x\right).
	\end{aligned}
\end{equation}
The size of the domain, described by $(x_l,y_l)$, is visualized in~\Cref{fig:domain_and_wound} as well. 

\begin{figure}[t]
    \centering
    \resizebox{1\textwidth}{!}{
    \begin{tikzpicture}[scale=1, >=Stealth]

    \node[draw, rectangle, minimum width=1.3cm] (branch_inputs) {
        \begin{tabular}{c}
            $D_F$ \\ 
            $\chi_F$ \\ 
            $D_c$ \\ 
            $k_F$ \\ 
            $a_c^I$
        \end{tabular}
    };

    \node[draw, rectangle, fill=golden_yellow, text opacity=1, above = 0cm of branch_inputs, minimum width=1.3cm] (branch_inputs2) {
        \begin{tabular}{c}
            $y_{cut}$ \\ 
            $x_m$ \\ 
            $y_m$ \\ 
            $x_{cut}$
        \end{tabular}
    };

    \node[draw, rectangle, right=1cm of branch_inputs, minimum height=1cm, minimum width=2cm] (branch_net) {Branch Net};

    \node[draw, rectangle, right=1.2cm of branch_net, minimum width=1.2cm] (branch_outputs) {
        \begin{tabular}{c}
            $b_1$ \\ 
            $\vdots$ \\ 
            $b_p$ \\ 
            \hline 
            $b_{p+1}$ \\ 
            $\vdots$ \\
            $b_{2p}$
        \end{tabular}
    };

    \node[draw, rectangle, below=1cm of branch_inputs, minimum width=1.3cm] (trunk_inputs) {
        \begin{tabular}{c}
            $t$ \\ 
            $x$ \\ 
            $y$ 
        \end{tabular}
    };

    \node[draw, rectangle, fill=cognac, fill opacity=0.3, text opacity=1, below=0cm of trunk_inputs, minimum width=1.3cm] (trunk_inputs2) {
        \begin{tabular}{c}
            $y_{cut}$ \\ 
            $x_{m}$ \\ 
            $y_{m}$ \\ 
            $x_{cut}$
        \end{tabular}
    };

    \node[draw, rectangle, right=1cm of trunk_inputs, minimum height=1cm, minimum width=2cm] (trunk_net) {Trunk Net};

    \node[draw, rectangle, right=1.2cm of trunk_net, minimum width=1.3cm, minimum height=1cm] (trunk_outputs) {
        \begin{tabular}{c}
            $c_1$ \\ 
            $\vdots$ \\ 
            $c_p$
        \end{tabular}
    };

    \node[circle, draw, right=1.3cm of branch_outputs] (inner_product1) {$\bigotimes$};
    \node[circle, draw, right=1.3cm of trunk_outputs] (inner_product2) {$\bigotimes$};

    \node[draw, rectangle, right=1cm of inner_product1, minimum height=1cm] (output1) {$\hat{u}_1(t,x,y)$};
    \node[draw, rectangle, right=1cm of inner_product2, minimum height=1cm] (output2) {$\hat{u}_2(t,x,y)$};

    \node[circle, draw, right=1cm of output1] (multiplication1) {$\times$};
    \node[circle, draw, right=1cm of output2] (multiplication2) {$\times$};

    \node[draw, rectangle, fill=sage, fill opacity=0.4, text opacity=1, above=1cm of multiplication1, minimum height=1cm] (sine1) {$\sin\left(\frac{\pi}{x_l}x\right)\cos\left(\frac{\pi}{2y_l}y\right)$};
    \node[draw, rectangle, fill=sage, fill opacity=0.4, text opacity=1, below=1cm of multiplication2, minimum height=1cm] (sine2) {$\sin\left(\frac{\pi}{y_l}y\right)\cos\left(\frac{\pi}{2x_l}x\right)$};    
    
    \node[draw, rectangle, right=1cm of multiplication1, minimum height=1cm] (final_output1) {$u_1(t,x,y)$};
    \node[draw, rectangle, right=1cm of multiplication2, minimum height=1cm] (final_output2) {$u_2(t,x,y)$};    
    
    \draw[->] (branch_inputs) -- (branch_net);
    \draw[->] (branch_net) -- (branch_outputs);

    \draw[->] (trunk_inputs) -- (trunk_net);
    \draw[->] (trunk_net) -- (trunk_outputs);

    \draw[->] ([yshift=1cm]branch_outputs.east) -- (inner_product1.west);
    \draw[->] ([yshift=-1cm]branch_outputs.east) -- (inner_product2.west);

    \draw[->] (trunk_outputs.east) -- (inner_product1.west);
    \draw[->] (trunk_outputs.east) -- (inner_product2.west);

    \draw[->] (inner_product1) -- (output1);
    \draw[->] (inner_product2) -- (output2);

    \draw[->] (output1) -- (multiplication1);
    \draw[->] (output2) -- (multiplication2);    

    \draw[->] (multiplication1) -- (final_output1);
    \draw[->] (multiplication2) -- (final_output2);

    \draw[->] (sine1.south) -- (multiplication1.north);
    \draw[->] (sine2.north) -- (multiplication2.south);    

    \node[above=1cm of branch_inputs2, xshift=-2cm] (proposed_model) {\textbf{Final model:}};
    \node[draw, rectangle, minimum width=0.5cm, minimum height=0.5cm, right=0.1cm of proposed_model] (white_box1) {};
    \node[right=0.05cm of white_box1] (plus1) {+};
    \node[draw, rectangle, fill=cognac, fill opacity=0.3, minimum width=0.5cm, minimum height=0.5cm, right=0.05cm of plus1] (cognac_box1) {};
    \node[right=0.05cm of cognac_box1] (plus2) {+};
    \node[draw, rectangle, fill=sage, fill opacity=0.4, minimum width=0.5cm, minimum height=0.5cm, right=0.05cm of plus2] (sage_box1) {};

    \node[right=1cm of sage_box1] (case1) {\textit{Case 1:}};
    \node[draw, rectangle, minimum width=0.5cm, minimum height=0.5cm, right=0.1cm of case1] (white_box2) {};

    \node[right=1cm of white_box2] (case2) {\textit{Case 2:}};   
    \node[draw, rectangle, minimum width=0.5cm, minimum height=0.5cm, right=0.1cm of case2] (white_box3) {};
    \node[right=0.05cm of white_box3] (plus3) {+};
    \node[draw, rectangle, fill=cognac, fill opacity=0.3, minimum width=0.5cm, minimum height=0.5cm, right=0.05cm of plus3] (cognac_box2) {};

    \node[right=1cm of cognac_box2] (case3) {\textit{Case 3:}};   
    \node[draw, rectangle, minimum width=0.5cm, minimum height=0.5cm, right=0.1cm of case3] (white_box4) {};
    \node[right=0.05cm of white_box4] (plus4) {+};
    \node[draw, rectangle, fill=golden_yellow, minimum width=0.5cm, minimum height=0.5cm, right=0.05cm of plus4] (yellow_box) {};

    \node[right=1cm of yellow_box] (case4) {\textit{Case 4:}};  
    \node[draw, rectangle, minimum width=0.5cm, minimum height=0.5cm, right=0.1cm of case4] (white_box5) {};
    \node[right=0.05cm of white_box5] (plus5) {+};
    \node[draw, rectangle, fill=sage, fill opacity=0.4, minimum width=0.5cm, minimum height=0.5cm, right=0.05cm of plus5] (sage_box2) {};
    \end{tikzpicture} }
    \caption{The DeepONet architecture used for predicting dermal displacement. The legend at the top describes four other cases against which we have tested our proposed model; see \Cref{tab:performance}}
    \label{fig:architecture}
\end{figure}
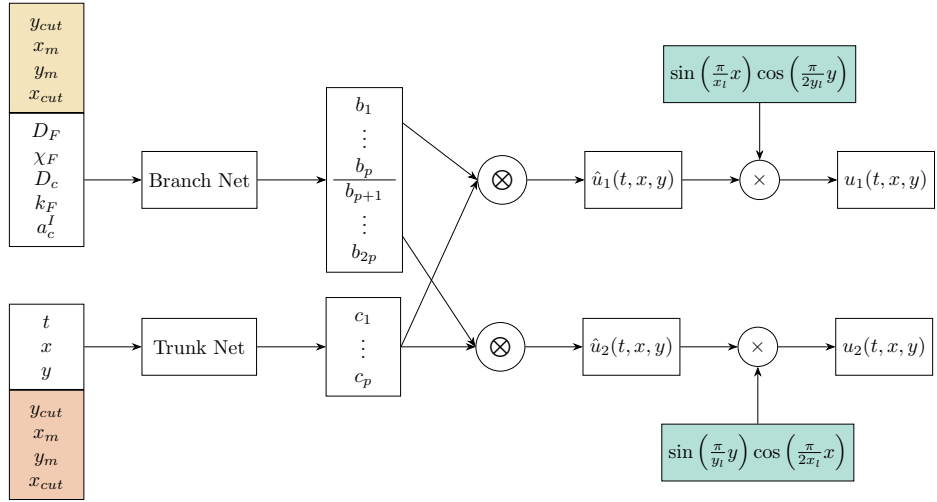

\section{Datasets and training} \label{sec:datasets_training}
\noindent
This section describes the approaches we used for generating datasets to train and test our model, as well as the specific training setup.

\subsection{Training set} \label{subsec:training}
\noindent
We generate a training set consisting of three different geometries for the initial wound shape: a rectangle (including a square), a rhombus (including a rotated square), and an ellipse (including a circle). We consider a quarter of the complete domain, assuming symmetry along the $x$- and $y$-axis, as shown in \Cref{fig:domain_and_wound}.
\par We run 750 finite element simulations, each run uniformly selecting one of the three initial wound shapes. For the wound size, we uniformly sample $x_{cut},y_{cut} \in (0,5)$ centimeters. The size of the complete domain is then determined as 
\begin{equation} \label{eq:x_l_y_l}
    \begin{aligned}
    x_l = 2.5 \cdot x_{cut}, \\
    y_l = 2.5 \cdot y_{cut},
    \end{aligned}    
\end{equation}
rounded to one decimal. For all simulations, we fix $t_{end} = 100$ days. The training data sampling strategy is as follows: for each finite element simulation, 10 time steps are uniformly sampled. Then, for each time step, 20 spatial coordinates from the domain are uniformly sampled. This gives us 200 coordinates $(t,x,y)$ per simulation, in which we record the displacements $(u_1, u_2)$. Each triplet $(t,x,y)$ is extended with the corresponding quadruple $(y_{cut}, x_m, y_m, x_{cut})$, which completes the input to the trunk network.  

\begin{table}[t]
\centering
\caption{Ranges for the values of the five parameters taken as input to the branch network} 
\begin{tabular}{cr@{\hspace{3pt}}c@{\hspace{3pt}}ll}
\toprule
\textbf{Parameter} & \multicolumn{3}{c}{\textbf{Range}}  & \textbf{Dimension}                            \\ \midrule
$D_F$              & $7.6167\cdot 10^{-7}$ & -- & $1.2 \cdot 10^{-6}$ & cm$^5$/(cells day) \\
$\chi_F$           & $(2$ & -- & $3)\cdot 10^{-3}$                      & cm$^5$/(g day)     \\
$D_c$              & $(2.22$ & -- & $3.2) \cdot 10^{-3}$              & cm$^2$/day          \\
$k_F$              & $8\cdot 10^6$ & -- & $1.08 \cdot 10^7$           & cm$^3$/(g day) \\
$a_c^I$            & $(0.9$ & -- & $1.1) \cdot 10^{-8}$               & g/cm$^3$              \\ \bottomrule
\end{tabular}
\label{tab:parameter_ranges}
\end{table}

\par For each numerical simulation, the five parameter values taken as input to the branch are sampled from uniform distributions with pre-specified ranges, as listed in \Cref{tab:parameter_ranges}. Furthermore, the values $(x_l, y_l)$, as given in~\Cref{eq:x_l_y_l}, are tabulated, as they are necessary for defining the sine augmentation; cf.~\Cref{eq:sine}. 

Following the procedure described above, each finite element simulation contributes the following quadruple to the training data:  

\begin{equation*}
    \begin{bmatrix}
        \begin{bmatrix}
            D_F, \chi_F, D_c, k_F, a_c \\
            D_F, \chi_F, D_c, k_F, a_c \\
            \vdots \\
            D_F, \chi_F, D_c, k_F, a_c \\
        \end{bmatrix}, 
        \begin{bmatrix}
            t_1, x_1,y_1,y_{cut},x_m, y_m, x_{cut} \\
            \vdots \\
            t_1, x_{20},y_{20},y_{cut},x_m, y_m, x_{cut} \\
            \vdots \\
            t_{10}, \hat{x}_1,\hat{y}_1,y_{cut},x_m, y_m, x_{cut} \\
            \vdots \\
            t_{10}, \hat{x}_{20},\hat{y}_{20},y_{cut},x_m, y_m, x_{cut} \\
        \end{bmatrix},
        \begin{bmatrix}
            u_1(t_1,x_1,y_1), u_2(t_1,x_1,y_1) \\
            \vdots \\
            u_1(t_1,x_{20},y_{20}), u_2(t_1,x_{20},y_{20}) \\
            \vdots \\
            u_1(t_{10},\hat{x}_1,\hat{y}_1), u_2(t_{10},\hat{x}_1,\hat{y}_1) \\
            \vdots \\
            u_1(t_{10},\hat{x}_{20},\hat{y}_{20}), u_2(t_{10},\hat{x}_{20},\hat{y}_{20}) \\
        \end{bmatrix},
        \begin{bmatrix}
            x_l, y_l \\
            x_l, y_l \\
            \vdots \\
            x_l, y_l \\
        \end{bmatrix}
    \end{bmatrix}. 
\end{equation*}
The first column is input to the branch, the second column is input to the trunk, the third column represents the target values, and the final column evaluates the sines and cosines in the sine augmentation step. As we perform 750 such numerical simulations, we obtain a training set containing $150\,000$ data points. We use an $80\,\%-20\,\%$ split for training and validation, respectively.

\subsection{Training setup}
\begin{figure}[t]
    \centering
    \includegraphics[width=0.5\columnwidth]{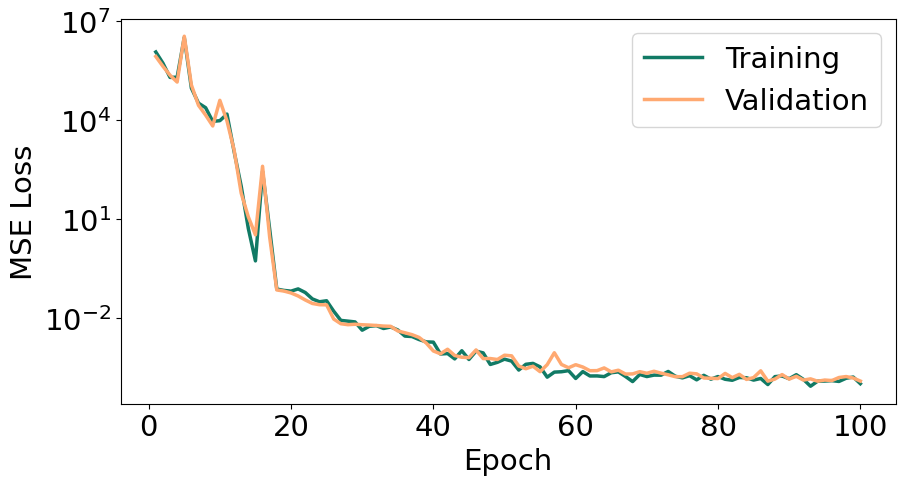}
    \caption{Training and validation losses of the proposed DeepONet. Dataset is based on 750 finite element simulations, using an 80\%-20\% split for training and validation, respectively}
    \label{fig:training_and_validation_loss}
\end{figure}

\noindent
Hyperparameter tuning was initially performed on a simplified model considering only one wound shape, and the resulting parameters are applied in this case as well (for more details, see the master's thesis \cite{husanovic2024}). For training our model, we use the Adam algorithm with standard backpropagation to minimize the MSE loss. The learning rate is set to 0.001. Through empirical testing, we found that a batch size of 100 yields satisfactory results. We train the network for 100 epochs, after which we observe that the loss stagnates around $10^{-3}$. Training the model for more epochs does not further decrease the loss. \Cref{fig:training_and_validation_loss} shows a plot of the training and validation losses.    

\subsection{Convex test set}
\noindent
We test the performance of our model on what we call the \textit{convex test set}. This is a more generic data set, containing convex combinations of the three basic shapes used for training. Let $s \in [a,b]$ and $\mathbf{r}_i = \mathbf{r}_i(s)$ represent a parametrization of the (piecewise) smooth boundary of basic wound shape $i$ ($i = 1$ corresponding to a rectangle, $i = 2$ to a rhombus, and $i=3$ to an ellipse), on which the DeepONet is trained. Then, a convex combination of the three parametrizations is given by
\begin{align} \label{eq:convex}
{\bf r}(s,\alpha_1,\alpha_2, \alpha_3) = \sum_{p=1}^3 {\bf r}_p(s) \, \alpha_p, 
\end{align} 
where $\alpha_p \ge 0$ and $\sum_{p=1}^3 \alpha_p = 1$. This parametrization provides the boundary of a new wound shape. \Cref{fig:sine_and_bnd} (b) gives a visualization. 
\par The convex test set it created by executing 150 finite element simulations. For each simulation, we create a convex combination of the three basic shapes. To this end, we first fix the sizes of the basic shapes by uniformly choosing $x_{cut},y_{cut} \in (0,5)$ centimeters. Then, we randomly select three weights that sum up to one and perform the multiplication as given in \Cref{eq:convex}.        
\par The remainder of the data sampling strategy is similar to the one used for generating the training data described in \Cref{subsec:training}. The only difference now is that for each finite element run, we sample all available time steps, and for each time step, we sample all available spatial coordinates in the domain. This procedure results in a convex test set containing 18\,035\,821 data points. 
\par We note that the quadruples $(y_{cut}, x_m, y_m, x_{cut})$ no longer uniquely define the initial wound shapes for this dataset. This is because it contains convex combinations of the three basic shapes, which do not have a straightforward geometric definition that can be captured by these four values alone. As a result, the network does not know beforehand the exact initial shapes of the samples in the convex test set.   

\section{Results} \label{sec:results}
\noindent
This section presents the results. We begin by examining the performance of the trained DeepONet on the convex test set, followed by a comparison to four other architectural setups. Lastly, we assess predictions over an extended period of one year.   

\begin{figure}[ht]
    \centering
    \begin{subfigure}[t]{.32\textwidth}
        \centering
        \includegraphics[width=\textwidth]{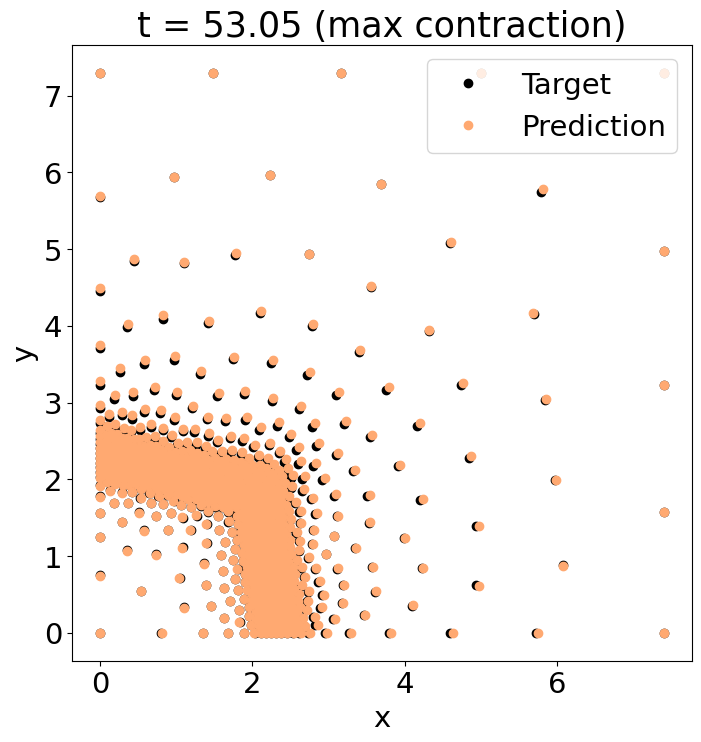}
        \caption*{(a)}
    \end{subfigure}
    \hfill
    \begin{subfigure}[t]{.32\textwidth}
        \centering
        \includegraphics[width=\textwidth]{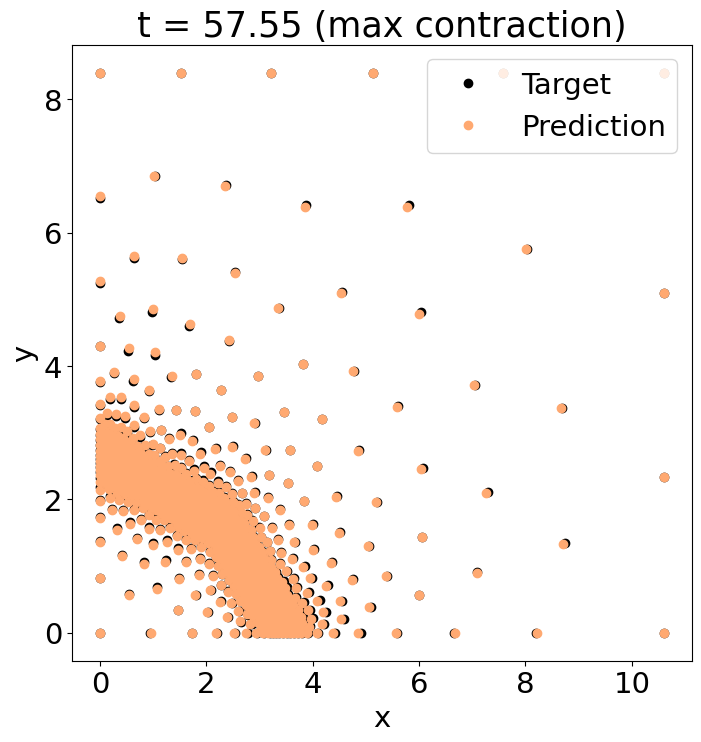}
        \caption*{(b)}
    \end{subfigure}    
    \hfill
    \begin{subfigure}[t]{.32\textwidth}
        \centering
        \includegraphics[width=\textwidth]{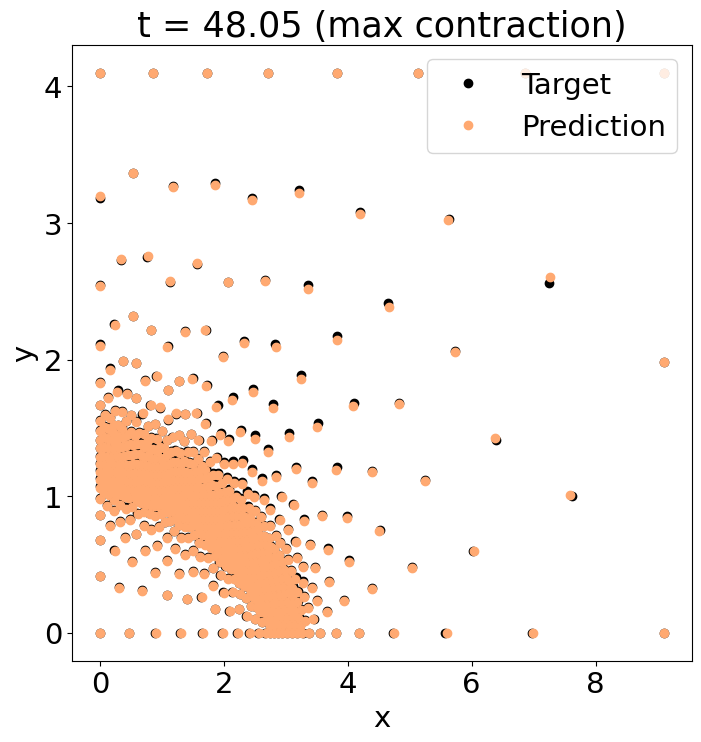}
        \caption*{(c)}
    \end{subfigure}
    \caption{Prediction vs. target for three samples from the convex test set, at their respective time of maximal contraction} 
    \label{fig:max_contraction_conv}
\end{figure}
      
\begin{figure}[ht]
    \centering
    \begin{subfigure}{.5\textwidth}
        \centering
        \includegraphics[width=0.7\textwidth]{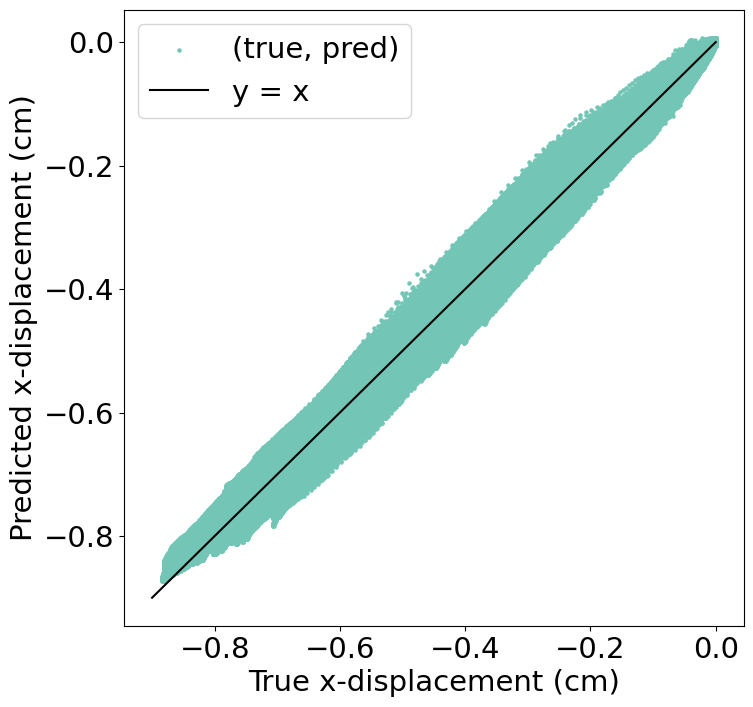}
        \caption*{(a)}
    \end{subfigure}%
    \begin{subfigure}{.5\textwidth}
        \centering
        \includegraphics[width=0.7\textwidth]{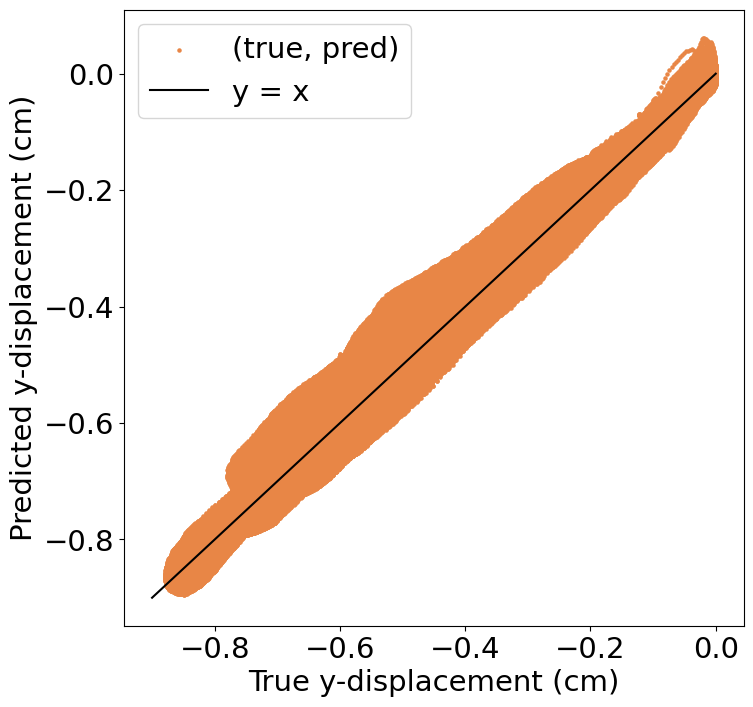}
        \caption*{(b)}
    \end{subfigure}
    \caption{True vs. predicted displacement in $x$- (a) and $y$-direction (b) on the convex test set}
    \label{fig:true_vs_predicted_conv}
\end{figure}

\begin{figure}[ht]
    \centering
    \begin{subfigure}[t]{.32\textwidth}
        \centering
        \includegraphics[width=\textwidth]{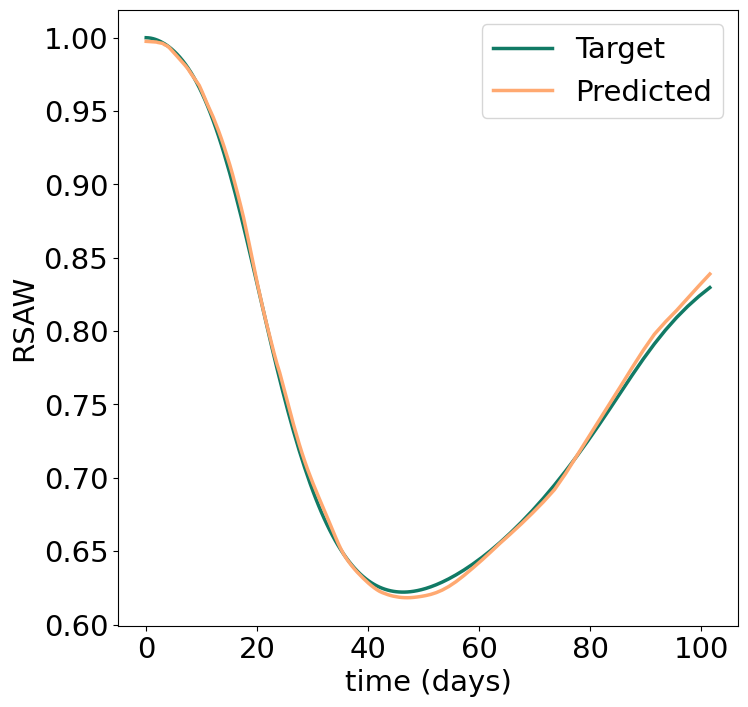}
        \caption*{(a)}
    \end{subfigure}
    \hfill
    \begin{subfigure}[t]{.32\textwidth}
        \centering
        \includegraphics[width=\textwidth]{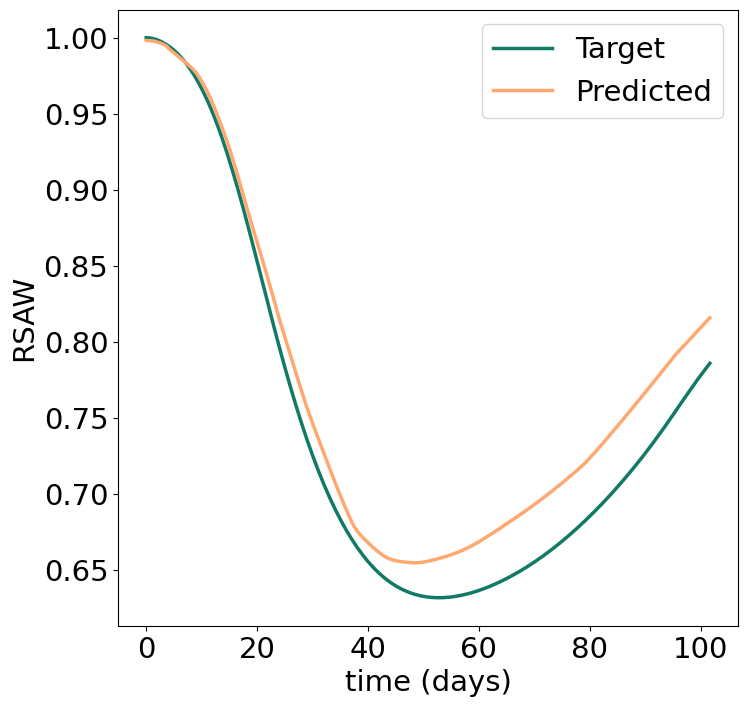}
        \caption*{(b)}
    \end{subfigure}    
    \hfill
    \begin{subfigure}[t]{.32\textwidth}
        \centering
        \includegraphics[width=0.98\textwidth]{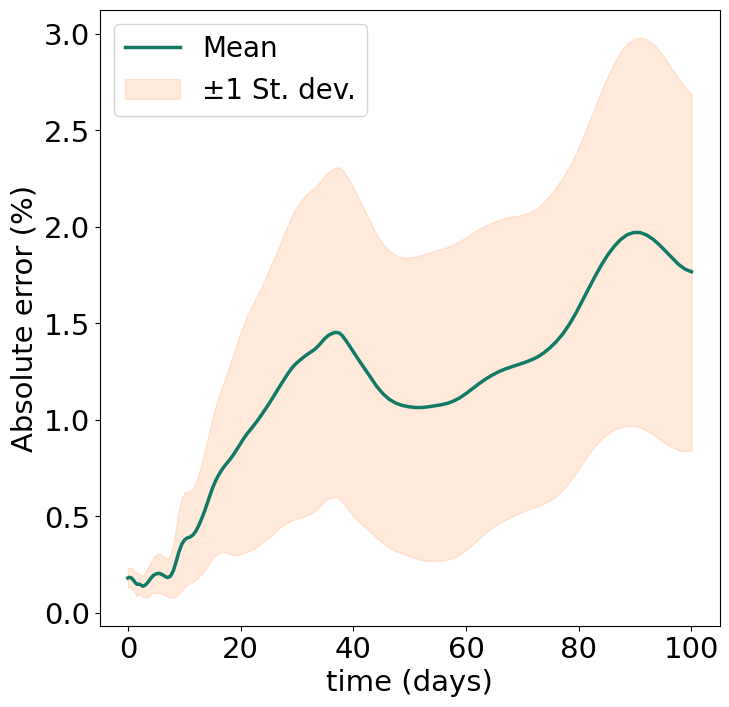}
        \caption*{(c)}
        \label{fig:AE}
    \end{subfigure}
    \caption{Best (a) and worst (b) prediction on the convex test set in terms of the RSAW. Mean and standard deviation of absolute error as a function of time (c). The mean is taken over all 150 samples in the convex test set} 
    \label{fig:rsaw_conv_AE}
\end{figure}

\subsection{Performance on convex test set}
\noindent
Our DeepONet predicts the $x$ and $y$ displacements in any given spatial coordinate in the domain, at any given time. We evaluate the trained model on each sample in the convex test set. \Cref{fig:max_contraction_conv} shows three different examples of the obtained predictions versus the corresponding targets, each at its respective time of maximal contraction. Note that we depict the actual coordinates after applying the displacements here, rather than the displacements. We observe good prediction abilities, especially at the wound boundary. Given that most of the spatial coordinates in the training set are located at the wound boundaries, it is expected that the network performs best in these regions. Further away from the wound, the prediction accuracy decreases slightly. As we move towards the top and right boundary of the domain, the displacements go to zero. We observe this in our predictions as well, which was ensured by the sine augmentation step during training.
\par Furthermore, given that our model is trained on the three basic wound shapes, \Cref{fig:max_contraction_conv} gives an indication of good generalization abilities. This is further validated in the scatter plots in \Cref{fig:true_vs_predicted_conv}, where the true displacements are plotted against the predicted displacements. The line $y=x$ indicates perfect predictions of the neural network model. The scatter plots contain all convex test set samples at all times steps, that is, around 18 million data points. We observe a strong concentration around the lines $y=x$, with no unusual outliers. This indicates that our model gives good predictions and that it generalizes well to convex combinations of the shapes is was trained on. We note that the scatter plots for $x$- and $y$-displacement are not exactly the same, as no symmetry was enforced; consequently the performance also differs slightly. 
For both $x$- and $y$-displacement, we observe that the model is best at predicting very small magnitude displacements (close to zero) and the largest magnitude displacements (around $-0.8$ cm). The former occurs near the boundary of the domain at all time steps, and at all spatial coordinates when $t$ is small. The latter occurs at the wound boundary at maximal contraction, for the largest wound sizes present in the convex set.    
\par Additionally, we evaluate our trained model on the wound boundary. We consider the same convex test set, but only the spatial coordinates on the wound boundary are taken. This enables us to determine the RSAW at each time step. Figs. \ref{fig:rsaw_conv_AE} (a) and (b) show the best and worst predictions in terms of the RSAW on the convex set. In the best case, the RSAW curve mostly overlaps the target, with a very slight underestimation of the maximal contraction and overestimation of the final value at $t=100$ days. In the worst case, the model predicts an earlier maximal contraction (49 days vs. 53) and retraction, with an overestimation of around 4\% less contraction from $t=55$ days and onwards. 
\par \Cref{fig:rsaw_conv_AE} (c) gives more insight in the behavior of the error over time. Here, the mean and standard deviation of the absolute error over the convex test set are depicted. We observe that the absolute error increases over time. The model is most accurate in the first 18 days, where the error is up to 0.5\% on average. After this time, the error, and the uncertainty thereof, start to increase. This coincides with very steep contraction. A notable peak can be observed around day 38, after which there is a local minimum (of 1.2\% on average) around day 50. The latter coincides with the average time of maximal contraction, indicating that the model is better at predicting the minimum RSAW (maximal displacement). The highest absolute error, with the largest uncertainty, is around day 90. For all times, the mean plus standard deviation of the error never exceeds 3\%. Since the model predicts displacements within the rough range of $[-0.9, 0]$ centimeters, a 3\% error corresponds to deviations on the order of tenths of millimeters.   

\subsection{Comparison to different architectures} 
\noindent
\Cref{tab:performance} presents the performance metrics of the DeepONet, where a comparison is made with four other architecture setups. These are variations of our proposed model (see \Cref{fig:architecture}). First, we define the performance measures used to determine the metrics in \Cref{tab:performance}. Subsequently, we give an interpretation of the results.      

\paragraph{Performance measures}
We include the goodness of fit ($R^2$) statistic, the average relative root mean squared error (aRRMSE), and the average relative error (aRelErr):

\begin{align}
    R^2 &= \frac{1}{2}\sum_{t=1}^2 \left[1 - \frac{\sum_{i=1}^n \left(u_t^{(i)}-\hat{u}_t^{(i)}\right)^2}{\sum_{i=1}^n \left(u_t^{(i)}-\bar{u}_t\right)^2}\right], \label{eq:R2} \\
    aRRMSE &= \frac{1}{2}\sum_{t=1}^2 \sqrt{\frac{\sum_{i=1}^n \left(u_t^{(i)}-\hat{u}_t^{(i)}\right)^2}{\sum_{i=1}^n \left(u_t^{(i)}-\bar{u}_t\right)^2}} \label{eq:aRRMSE}, \\
    aRelErr &= \label{eq:aRelErr} 
    \frac{1}{2} \sum_{t=1}^{2} \frac{1}{n_t} \sum_{i=1}^{n} \mathbb{I}(u_t^{(i)} \neq 0) \cdot \left| \frac{u_t^{(i)} - \hat{u}_t^{(i)}}{u_t^{(i)}} \right|,
\end{align}
where 
\begin{itemize}
    \item $u_t^{(i)}$ is the true displacement (rounded to one decimal place in \Cref{eq:aRelErr}),
    \item $\hat{u}_t^{(i)}$ is the predicted displacement (rounded to one decimal place in \Cref{eq:aRelErr}),
    \item $\bar{u}_t$ is the average true displacement,
    \item $\mathbb{I}(u_t^{(i)} \neq 0)$ is the indicator function, which is 1 if $u_t^{(i)} \neq 0$ and zero otherwise,  
    \item $n_t$ is the number of samples where $u_t^{(i)} \neq 0$,
    \item $n$ is the total number of samples.
\end{itemize}

\begin{table}[t]
\caption{Performance of the proposed DeepONet compared to four other setups (see \Cref{fig:architecture} for a visualization). Boldface indicates the best performance}
\resizebox{\textwidth}{!}{%
\begin{tabular}{@{}lccccc@{}}
\toprule
\textbf{\begin{tabular}[c]{@{}l@{}}Performance\\ metric\end{tabular}} & \begin{tabular}[c]{@{}c@{}}\textbf{No sine aug. and} \\ \textbf{no initial shape info} \\ (Case 1)\end{tabular} & \begin{tabular}[c]{@{}c@{}}\textbf{Initial shape info}\\ \textbf{to branch} \\ (Case 2) \end{tabular} & \begin{tabular}[c]{@{}c@{}}\textbf{Initial shape info} \\ \textbf{to trunk} \\ (Case 3)\end{tabular} & \begin{tabular}[c]{@{}c@{}}\textbf{Sine} \\ \textbf{augmentation} \\ (Case 4)\end{tabular} & \begin{tabular}[c]{@{}c@{}}\textbf{Sine aug. and initial} \\  \textbf{shape info to trunk} \\ (Proposed model)\end{tabular} \\ \midrule
$R^2$                                                                 & $0.90040$                                                                            & $0.89613$                                                                    & $0.99155$                                                                   & $0.89083$                                                             & $\mathbf{0.99437}$                                                                               \\
aRRMSE                                                                & $0.31560$                                                                            & $0.32229$                                                                    & $0.09166$                                                                   & $0.33032$                                                             & $\mathbf{0.07497}$                                                                               \\
aRelErr                                                               & $0.12160$                                                                            & $0.11673$                                                                    & $0.04432$                                                                   & $0.15050$                                                             & $\mathbf{0.03429}$                                                                               \\ \bottomrule
\end{tabular}}
\label{tab:performance}
\end{table}

\paragraph{Interpretation}
A value of $R^2 = 1$ indicates perfect predictions, whereas $R^2 = 0$ indicates that the model is not better than always returning the expected values. For our proposed DeepONet, we find $R^2 = 0.99$. This shows that the model can accurately predict dermal displacement. Furthermore, we find an aRRMSE of $0.08$, showing excellent performance according to Despotovic
et al. \cite{despotovic2016evaluation}. Lastly, the aRelErr is found to be only $0.03$, which further indicates that the model can excellently reproduce the finite element simulations. 
\par We compare the proposed DeepONet architecture to four variations thereof. \Cref{fig:architecture} shows which blocks are included in the different cases, corresponding to the columns of \Cref{tab:performance}. We observe that, compared to all variations, the proposed model gives the best performance. The "skeleton" DeepONet, where the model does not have any information about the initial wound shape and no sine augmentation is applied (case 1), performs significantly worse. We find $R^2 = 0.90$, $\textrm{aRRMSE} = 0.32$, and $\textrm{aRelErr} = 0.12$. Adding initial wound shape info in the form of the quadruple $(y_{cut}, x_m, y_m, x_{cut})$ to the branch (case 2), even gives slightly worse performance. The largest gain is when this initial shape info is added to the trunk (case 3). Here we see a decrease of $71 \%$ in the aRRMSE and a decrease of $63.6\%$ in the aRelErr, compared to case 1. This suggests that initial shapes should be input to the trunk, as this component of the network is responsible for learning the basis functions. Improved representation of basis functions often leads to a more accurate approximation of the solution. It is also the most logical choice because the initial shape info is more closely related to the spatial coordinates, which are also fed into the trunk, rather than to the parameters, which are input to the branch. Case 4, where the skeleton DeepONet only has the addition of the sine augmentation blocks, is the worst performing one with $R^2 = 0.89$, $\textrm{aRRMSE}=0.33$, and $\textrm{aRelErr} = 0.15$. From this, we conclude that the network cannot sufficiently learn the dermal displacements without knowledge of the initial shape of the wound. Additionally, sine augmentation provides an improvement in performance only when the model already performs sufficiently without it. Indeed, if we combine cases 3 and 4, we see that our proposed model overall results in the best performance. 

\subsection{Prediction after $t=100$}
\begin{figure}[t]
    \centering
    \begin{subfigure}{.5\textwidth}
        \centering
        \includegraphics[width=0.7\textwidth]{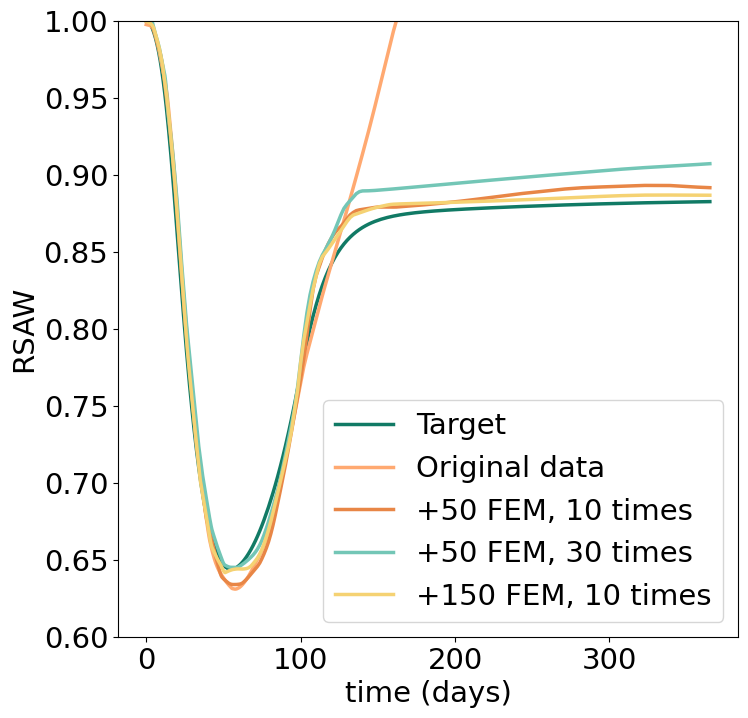}
        \caption*{(a)}
        \label{fig:rsaw_comparison_best}
    \end{subfigure}%
    \begin{subfigure}{.5\textwidth}
        \centering
        \includegraphics[width=0.7\textwidth]{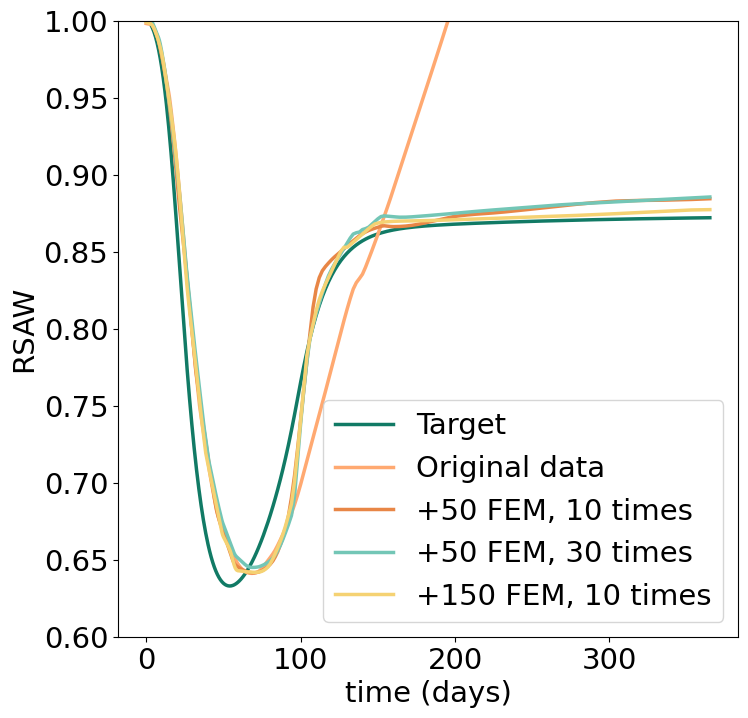}
        \caption*{(b)}
        \label{fig:rsaw_comparison_worst}
    \end{subfigure}
    \caption{Best (a) and worst (b) prediction in terms of the RSAW on a small convex test set, comparing different training data sampling strategies. The legend denotes the number of finite element simulations added to the training set, and the number of time steps sampled per simulation}
    \label{fig:best_worst_rsaw_sampling}
\end{figure}

\noindent
So far, we have taken $t\in[0,100]$ days for all predictions to save time when generating data using numerical simulations. 
However, we are interested to see if our final DeepONet is also capable of predicting $t\in[0,365]$ days. \Cref{fig:best_worst_rsaw_sampling} shows two results of evaluating the DeepONet trained on $t \in [0,100]$ days (the original data), on an entire year. We observe that the predicted RSAW is linearly increasing after $t = 100$ days, and no asymptotic value is reached, which means the network predicts a scar that continues to grow in size. This indicates the need for adding data to the training set with $t_{end} = 365$. We experiment with different data sampling strategies. Specifically, we compare the effects of incorporating data from additional finite element simulations, to drawing more data points from existing samples. 
Practically, this comparison reflects the difference between including data on more burn patients versus collecting more data per patient. 
\par We generate a new convex test set for evaluating the model's performance, consisting of 50 finite element runs with $t_{end} = 365$. \Cref{fig:best_worst_rsaw_sampling} shows the results from three different training data sampling strategies. For each case, we initialize the network with the learned parameters from the final DeepONet trained on $t\in [0,100]$, as a good initial guess to speed up the training. In this way, the network embeds prior knowledge, instead of having to start learning from scratch.
\par We compare the scenario where 50 finite element simulations are added to the training set, with 30 time steps sampled for each (scenario 1), versus the scenario where 150 runs are added, with 10 time steps sampled (scenario 2). Note that both result in the same number of additional data points (30\,000). In both the best and worst cases, we observe that scenario 2 performs better. Sampling more data per existing finite element run does not lead to much improvement. In fact, in~\Cref{fig:best_worst_rsaw_sampling} (a), we even observe a decrease compared to selecting 10 time steps. 
Adding new, unique data results in a performance increase. In the best case, the predicted RSAW at $t=365$ is now only 0.47\,\% higher than the target. In the worst case, this is 0.6\,\%. 
\par The above results show that adding more runs is preferred over adding more data points within a run, possibly because adding runs increases the diversity of the training set. It exposes the network to a wider variety of features, which can enable it to generalize better. On the other hand, introducing more data per existing samples can introduce redundancy, as the network is repeatedly exposed to the same patterns. Since these patterns were already learned in the original data, additional exposure does not result in significant improvement and can lead to overfitting. This means that the network becomes too specialized in recognizing the patterns in the training data, reducing its ability to generalize to new, unseen data.  

\subsection{Speedup}
\noindent In this study, all computations were performed on the DelftBlue supercomputer, provided by the Delft High Performance Computing Centre \cite{DHPC2024}. 
For the finite element simulations, we used a single node equipped with an Intel Xeon Gold 6248R CPU, using one core and 2 GB of memory. 
Executing the finite element-based model for a one-year prediction required, on average, 4.7 minutes of computing time per sample. In contrast, using the same computational resources, the DeepONet predicted the displacement field in an average of 2.2 seconds per sample, resulting in a speedup of about 128 times. Furthermore, the DeepONet has the advantage that it can be evaluated on a GPU. Using one NVIDIA A100 GPU on DelftBlue, the average prediction time per sample decreased to 1.2 seconds, which gives an additional speedup of almost two times.   
This shows that, in addition to being effective, DeepONets can serve as fast surrogates for traditional finite element methods in simulating post-burn contraction.
\par However, previous studies \cite{Egberts2022, Egberts2023} reported larger speedups compared to the numerical model. In \cite{Egberts2023}, MLPs were trained for the much simpler task of predicting the RSAW and wound boundary over time, given 25 parameters from the 2D morphoelastic model. Here, a total of 21 spatial points were used to describe a fixed initial wound shape. Since this differs from our approach of predicting the entire displacement field over time for multiple wound shapes, even shapes that were not 'seen' in the training data, a direct one-to-one comparison is not possible. The CPU time needed for our DeepONet to evaluate only the nodes on the wound boundary (100 in total) for each $t\in[0,365]$ is, on average, 0.27 seconds per sample. If we compare this to the execution time of the numerical model, we find a speedup of about 1044 times. 
\par The 2.2 seconds reported above for the DeepONet corresponds to evaluating the DeepONet on all nodes $(t,x,y)$ in the space-time grid. This means that, per sample evaluation, many forward passes need to be made until the network has processed all grid points in the domain, for each $t\in [0,365]$. For example, one sample from the convex test used for one-year predictions contains an average of $180\,010$ points $(t,x,y)$. Consequently, with a batch size of $50$, the model has to be evaluated approximately $3600$ times per sample.  
This is different from the previous approaches, where the MLP had to be evaluated only once to predict all quantities for each case.
However, the flexibility regarding geometry realized by our DeepONet is a significant advancement, while still achieving a speedup 
of about 128 on CPU and 235 on GPU. 

\section{Conclusions and discussion} \label{sec:conclusion}
\noindent
We trained a DeepONet on three different initial wound shapes to accurately predict the dermal displacement across the entire domain over time. The DeepONet consists of a branch and trunk network. The branch network takes as input five parameters from the morphoelastic model 
describing burn injuries: the cell diffusion constant, the chemotactic constant, the signaling diffusion constant, the cell differentiation rate, and the half-maximum cell division enhancement rate. The trunk takes as input a coordinate $(t,x,y)$ where we want to evaluate the solution, along with information about the initial wound shape and size given by the quadruple  $(y_{cut}, x_m, y_m, x_{cut})$. The final DeepONet output is the result of a dot product of the branch and trunk network outputs, such that the five input parameters are mapped to the displacement fields for 100 post-burn days.
\par We evaluated the performance of the trained DeepONet on a test set including finite element simulations based on convex combinations of the three basic wound shapes. We found that the predictions closely matched targets, demonstrating good generalization. The scatter plots in \Cref{fig:true_vs_predicted_conv} confirm this, showing a strong concentration around $y=x$ with minimal outliers. The model slightly underperformed on mid-range $y$-displacements, but excelled in predicting both small and large displacements. The error analysis in \Cref{fig:AE} showed a gradual increase in absolute error over time, peaking around day 90. However, the mean plus standard deviation of the error never exceeds 3\%, which corresponds to deviations on the order of tenths of millimeters. From a practical medical perspective, these are very accurate predictions.  
\par We compared our various DeepONet architectures on the convex test set, cf.~\Cref{fig:architecture}. The final model, incorporating initial wound shape info and sine augmentation, achieved the best performance with an $R^2$ score of $0.99$, an aRRMSE of $0.08$, and an aRelErr of $0.03$. This confirms the model's strong generalization abilities. The results support that initial shape info should be input to the trunk, and sine augmentation enhances performance when the model is already well-tuned. 
\par Finally, we extended the DeepONet predictions to one year, rather than 100 days. We expanded the training set with 50 new finite element simulations, where $t_{end}=365$. We compared the RSAW curves for two scenarios: one where 50 finite element simulations were added to the training set, with 30 time steps sampled in $t \in (100, 365]$, and another where 150 runs were added, with 10 time steps each. We found that adding new, unique data (the first scenario) results in the greatest improvement. Practically, these findings suggest that including more patients with less data per patient is preferable to having fewer patients with more data per patient. This is likely because introducing new information provides a more varied learning experience, whereas additional data points from existing samples tend to add redundancy, offering less benefit. 
\par We conclude that we have successfully trained a DeepONet as a finite element surrogate, capable of accurately predicting the post-burn dermal displacement field. The training encompassed multiple initial wound shapes, with the network demonstrating strong generalization on convex combinations of these basic shapes. Furthermore, by adding a limited number of samples from one-year finite element simulations into the training set, the network achieved reasonable predictions of the RSAW over the entire year. We anticipate that including more data with $t_{end} = 365$ will further enhance prediction accuracy. 
\par As a final conclusion, we experimentally found that the speedup of a simulation run using the DeepONet is about a factor 128 relative to the original finite element simulations (235 when evaluated on a GPU). While this speedup is lower than that achieved in~%
\cite{Egberts2022, Egberts2023} for MLP models on a fixed geometry,
the proposed DeepONet allows for geometrical flexibility that was not possible with the MLPs employed in previous works. 
\par A key aspect of our research was the novelty of incorporating multiple initial wound shapes. A promising direction for further research is to investigate the application of more complex and realistic initial wound geometries, where, for example, no symmetry in $x$ and $y$ is present. One possible approach is to include convolutional layers before the trunk network, allowing for the direct input of an image of the initial wound shape, thereby increasing the model applicability. Additionally, an alternative to sine augmentation is to penalize discrepancies between the predicted and true displacements on the domain boundary within the loss function, adopting a PINN-like approach for softer enforcement of the boundary conditions \cite{raissi_physics-informed_2019}. Further integration of DeepONets with PINNs can also ensure that all predicted displacements adhere more closely to the governing PDEs, which now is not necessarily the case \cite{wang2021learning}. This could solve potential discrepancies between the predictions and the underlying physical laws, leading to more reliable outcomes. \par Another direction for further research is to including a broader range of patient-specific parameters as input to the branch. This approach could provide a more comprehensive and applicable model; however, it would likely introduce larger variability in the solutions, since some parameters (e.g., collagen) are particularly sensitive. This would present an interesting opportunity to test whether the model can manage this increased complexity. Consequently, the architecture may need to be extended or improved. 

\section*{Datasets and code}
\noindent All relevant code is available at \href{https://github.com/Selma24/DeepONet-contraction}{github.com/Selma24/DeepONet-contraction}. All relevant data is available on the 4TU.Centre for Research Data at \href{https://doi.org/10.4121/69d1aefc-a01d-4280-8b32-5c8420d9a2a3}{doi.org/10.4121/69d1aefc-a01d-4280-8b32-5c8420d9a2a3}. 

\section*{Acknowledgements}
\noindent The authors acknowledge the use of computational resources of the DelftBlue supercomputer, provided by Delft High Performance Computing Centre (https://www.tudelft.nl/dhpc) \cite{DHPC2024}. The authors, Egberts and Vermolen, are grateful for the financial support from the Dutch Burns Foundation
under projects 17.105, 22.104 and PPS 22.01.

\section*{Declaration}
\noindent The authors declare that they have no conflict of interest.

\bibliographystyle{plain}
\bibliography{bibliography}

\newpage
\appendix
\section*{Appendix}
\noindent \Cref{tab:Parameters_stab_1d} lists the fixed parameter values used for the finite element simulations that were not taken as input to the DeepONet's branch network. 
\begin{table}[h]
    \centering
    \caption{Overview of the fixed parameter values used for the numerical simulations. 'NC' indicates that the parameter value is a consequence of the chosen values for other parameters.}  
    \label{tab:Parameters_stab_1d}  
    \begin{tabular}{@{}lllc@{}}  
        \toprule
        \textbf{Symbol} & \textbf{Value} & \textbf{Dimension} & \textbf{Reference} \\
        \midrule  
        $k_c$ & $4\cdot10^{-13}$ & g/(cells day) & \cite{Olsen1995} \\
        $r_F$ & $9.24\cdot10^{-1}$ & cm$^{3q}$/(cells$^q$ day) & \cite{Alberts1989,Gosh} \\
        $r^{\text{max}}_F$ & $2$ & - & \cite{Strutz2001} \\
        $k_\rho$ & $7.6\cdot10^{-8}$ & g/(cells day) & [NC] \\
        $k_\rho^{\text{max}}$ & $10$ & - & \cite{Olsen1995} \\
        $a_c^{II}$ & $10^{-8}$ & g/cm$^3$ & \cite{Olsen1995} \\
        $a_c^{III}$ & $2\cdot10^8$ & cm$^3$/g & \cite{Overall1991} \\
        $a_c^{IV}$ & $10^{-9}$ & g/cm$^3$ & \cite{Roberts1986} \\
        $\eta^I$ & $2$ & - & \cite{Rudolph1991} \\
        $\eta^{II}$ & $5\cdot10^{-1}$ & - & \cite{Koppenol2017a} \\
        $\kappa_F$ & $10^{-6}$ & cm$^3$/cells & \cite{VandeBerg1989} \\
        $q$ & $-4.151\cdot10^{-1}$ & - & [NC] \\
        $\delta_N$ & $2\cdot10^{-2}$ & /day & \cite{Olsen1995} \\
        $\delta_M$ & $6\cdot10^{-2}$ & /day & \cite{koppenol2017mathematical} \\
        $\delta_c$ & $5\cdot10^{-4}$ & cm$^6$/(cells g day) & \cite{Olsen1995} \\
        $\delta_\rho$ & $6\cdot10^{-6}$ & cm$^6$/(cells g day) & \cite{koppenol2017mathematical} \\
        $\overline{N}$ & $10^{4}$ & cells/cm$^3$ & \cite{Olsen1995} \\
        $\overline{M}$ & $0$ & cells/cm$^3$ & \cite{Olsen1995} \\
        $\overline{c}$ & $0$ & g/cm$^3$ & \cite{koppenol2017mathematical} \\
        $\overline{\rho}$ & $1.125\cdot10^{-1}$ & g/cm$^3$ & \cite{Olsen1995} \\
        $\rho_t$ & 1.09 & g/cm$^3$ & \cite{Icrp2010} \\
        $\mu_1$ & $10^2$ & (N day)/cm$^2$ & \cite{Koppenol2017a} \\
        $\mu_2$ & $10^2$ & (N day)/cm$^2$ & \cite{Koppenol2017a} \\
        $E$ & $32$ & N/((g cm)$^{0.5}$) & \cite{Egberts_2021s} \\
        $\xi$ & $5\cdot10^{-2}$ & (N g)/(cells cm$^2$) & \cite{Maskarinec, Wrobel2002} \\
        $R$ & $9.95\cdot10^{-1}$ & g/cm$^3$ & \cite{Koppenol2017a} \\
        $\zeta$ & $4\cdot10^2$ & cm$^6$/(cells g day) & \cite{Koppenol2017a} \\
        $\nu$ & $4.9\cdot10^{-1}$ & - & \cite{Egberts_2021s} \\ 
        $\widetilde{N}$ & $2\cdot10^{3}$ & cells/cm$^3$ & [NC] \\
        $\widetilde{c}$ & $10^{-8}$ & g/cm$^3$ & \cite{Egberts_2021s} \\
        $\widetilde{\rho}$ & $1.13\cdot10^{-2}$ & g/cm$^3$ & [NC] \\
        \bottomrule
    \end{tabular}
\end{table}

\end{document}